\documentclass{article}

\usepackage{neurips_2024}

\makeatletter
\@submissionfalse  
\makeatother

\nolinenumbers

\usepackage[utf8]{inputenc} 
\usepackage[T1]{fontenc}    
\usepackage{hyperref}       
\usepackage{url}            
\usepackage{booktabs}       
\usepackage{amsfonts}       
\usepackage{nicefrac}       
\usepackage{microtype}      
\usepackage{xcolor}         
\usepackage{graphicx}
\usepackage{float} 
\usepackage{subfigure}
\usepackage{wrapfig}
\usepackage{multirow}
\usepackage{algorithm, algorithmic}
\usepackage{hyperref}

\title{Learning to Select Cutting Planes in Mixed Integer Linear Programming Solving}

\author{%
  Xuefeng Zhang\textsuperscript{1} \\
  \And
  Liangyu Chen\textsuperscript{1} \\
  \And
  Zhengfeng Yang\textsuperscript{1}\\
  \And
  Zhenbing Zeng\textsuperscript{2} \\
}

\begin{document}

\maketitle

\begin{center}
  \textsuperscript{1}Software Engineering Institute, East China Normal University \\
  \textsuperscript{2}College of Sciences, Shanghai University \\
\end{center}

\begin{abstract}
Cutting planes (cuts) are crucial for solving Mixed Integer Linear Programming (MILP) problems. Advanced MILP solvers typically rely on manually designed heuristic algorithms for cut selection, which require much expert experience and cannot be generalized for different scales of MILP problems. Therefore, learning-based methods for cut selection are considered a promising direction. State-of-the-art learning-based methods formulate cut selection as a sequence-to-sequence problem, easily handled by sequence models. However, the existing sequence models need help with the following issues: (1) the model only captures cut information while neglecting the Linear Programming (LP) relaxation; (2) the sequence model utilizes positional information of the input sequence, which may influence cut selection.
To address these challenges, we design a novel learning model \emph{HGTSM} for better select cuts. We encode MILP problem state as a heterogeneous tripartite graph, utilizing heterogeneous graph networks to fully capture the underlying structure of MILP problems. Simultaneously, we propose a novel sequence model whose architecture is tailored to handle inputs in different orders. Experimental results demonstrate that our model outperforms heuristic methods and learning-based baselines on multiple challenging MILP datasets. Additionally, the model exhibits stability and the ability to generalize to different types of problems.
\end{abstract}

\section{Introduction}
Combinatorial optimization is a class of problems that find the optimal object in a limited set of objects. The methods and techniques related to combinatorial optimization have been widely developed in fields such as operations research,  computer science, etc., and have important applications in engineering, economics, and many other areas~\cite{bengio2021machine}. 
Among various combinatorial optimization problems, mixed-integer linear programming (MILP) is an extremely important one. 
Generally, $MILP=(\mathbf{A},\mathbf{b},\mathbf{c},\mathcal{I})$ is formally defined as:
\begin{equation}
z^{\ast}=\min\{\mathbf{c}^T\mathbf{x}|\mathbf{A}\mathbf{x}\le \mathbf{b},\mathbf{x}\in \mathbb{R}^n,x_j\in \mathbb{Z},\forall j\in \mathcal{I}\},
\end{equation}
where $\mathbf{A} \in \mathbb{R}^{m\times n}$ is referred to the constraint coefficient matrix, $\mathbf{b} \in \mathbb{R}^m$ is called the constraint right-hand-side vector, $\mathbf{c} \in \mathbb{R}^n$ is called the objective coefficient vector, and $\mathcal{I} \subseteq {1, 2, \ldots, n}$. Specifically, the decision variable vector $\mathbf{x}$ is required to be integers. Under the constraints mentioned above, the vector $\mathbf{x}^{\ast}$ that minimizes $\mathbf{c}^T\mathbf{x}$ to obtain the minimum value $z^{\ast}$ is the optimal solution. $\mathbf{X}_{MILP} = \{\mathbf{x} \in \mathbb{R}^n | \mathbf{A}\mathbf{x} \le \mathbf{b}, \mathbf{x} \in \mathbb{R}^n, x_j \in \mathbb{Z}, \forall j \in \mathcal{I}\}$ is called the feasible solution set, which refers to the set of vectors $\mathbf{x}$ that satisfy the constraints.

With the continuous expansion of the scale of MILP problems in the real world, how to solve MILP accurately and quickly has become very important. The typical algorithm for the exact solution of MILP, is the Branch and Bound (B\&B) method~\cite{land2010automatic}.
B\&B method uses the principle of divide and conquer, continuously divides the solution space to form a tree-like structure, and achieves a similar effect to enumeration. Once the tree is entirely explored, the exact solution to the problem can be obtained. In order to further improve the solving efficiency, solvers usually integrate the method of cutting planes (cuts)~\cite{gomory1960algorithm} into the B\&B method to form the Branch and Cut (B\&C) method~\cite{padberg1991branch}. The cut method can tighten the feasible domain before solving specific subproblems in the B\&B process while ensuring that no integer feasible solutions are affected, which undoubtedly enhances the solving efficiency. However, adding cuts to MILP arbitrarily may bring extra computational burdens and lead to numerical instability~\cite{wesselmann2012implementing,dey2018theoretical}. Therefore, how to add cuts properly has a significant impact on the solving efficiency of MILP. Unfortunately, although B\&C has become the universal framework in modern MILP solvers, the process of adding cuts still relies heavily on a large number of manually-designed heuristic methods~\cite{bengio2021machine,tang2020reinforcement}, including the generation of cuts, the number of cuts added, the selection of which cuts to add, and the order in which cuts are added, and so on. Heuristic methods often rely on specific problem characteristics or expert experience and cannot be applied to a wider range of problem-solving. Therefore, applying data-driven learning-based methods to the selection process of cuts instead of expensive and problem-specific heuristic methods is a meaningful and promising step for MILP solving.

Based on this idea, this paper focuses on the cut selection problem and proposes a novel model \emph{HGTSM}(Heterogeneous-Graph-Transformer Sequence Model) for learning the cut selection strategy. First, in order to fully capture the underlying information of MILP and enable the model to learn more general selection strategies, we describe MILP state as a tripartite graph and extract feature information using custom PySCIPOpt functions with the open-source solver SCIP~\cite{bestuzheva2021scip}. Second, building upon this, we utilize a heterogeneous graph neural network to process input information, 
where different types of nodes and edges are mapped to type-specific feature spaces, allowing for message passing and aggregation among nodes of different types. 
Third, we generate the selected sequence of cuts using a sequence model. We employ a reward-based learning approach to train the model with the expectation of learning a cut selection strategy, which is better than hand-designed heuristic methods. Finally, we conduct experiments on six public datasets, including two of medium difficulty and four of hard difficulty from MIPLIB 2017~\cite{gleixner2021miplib} and the NeurIPS 2021 ML4CO competition~\cite{ecoleai2021}. Experimental results demonstrate that the selected cut sequence by our model can better speed up the MILP-solving process, indicating the strong potential of learning-based models in improving solver performance. Furthermore, our model exhibits better generalization to different types of MILP problems compared to baseline methods and incurs acceptable performance loss, showcasing the prospects of seeking more general MILP-solving methods.

Specifically, our contributions can be summarized as follows:
\begin{itemize}
\item 
We recognize the heterogeneity of MILP denoted as a tripartite graph and choose heterogeneous graph networks to fully leverage this heterogeneity, thereby better capture the underlying information of MILP.
\item 
Based on extensive empirical results, we observe that the cuts selection methods based on sequence models captures positional information from input sequences, which is detrimental to problem solving. Therefore, we purposefully improve the sequence model to disregard positional information from input sequences, thereby enhance problem-solving stability.
\item 
We conduct experiments on six challenging datasets. The results indicate that our model not only outperforms the baselines in terms of performance but also significantly enhances problem-solving stability. Additionally, the model demonstrates good generalization.
\end{itemize}

The remainder of this paper is organized as follows. In Section~\ref{sec2}, we introduce the generation of motivation and demonstrate it with experiments. In Section~\ref{sec3}, we describe the encoding method of MILP states and the model architecture. In Section~\ref{sec4}, we provide a detailed explanation of the experimental setting and results. Finally, in Section~\ref{sec5}, we provide a summary and outlook.

\section{Motivation}
\label{sec2}
\begin{wrapfigure}{r}{0.5\textwidth}
    \vspace{-20pt}
    \centering
    \includegraphics[width=0.5\textwidth]{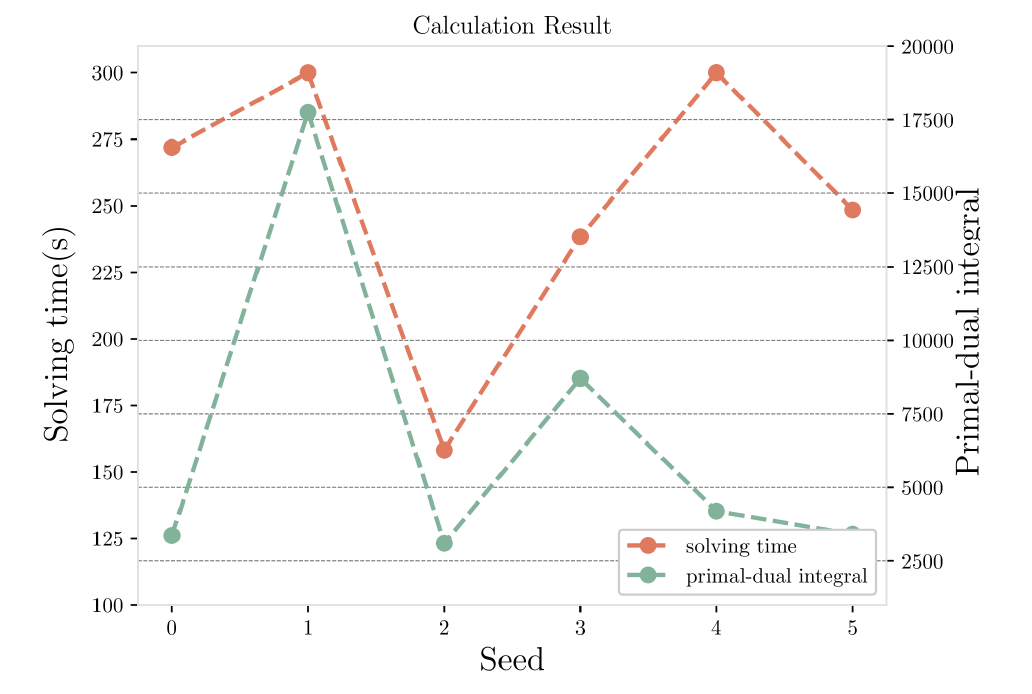}
    \vspace{-20pt}
    \caption{\footnotesize The time taken by the sequence model to solve the problem and the primal-dual gap integral(PD integral) when using input sequences shuffled with different random number seeds. }
    \label{motivation_f}
    \vspace{-15pt}
\end{wrapfigure}

The state-of-the-art (SOTA)  learning-based methods deemed cut selection as a sequence-to-sequence problem, handled with sequence models easily. However, the sequence models emphasize the potential positional information in input sequences, which empirically invokes a new issue: the different input orders of the candidate cuts have a big impact on solving performance. The intuitive idea is that for the same problem state and the same candidate cuts, the model should choose the same sequence of cuts, even if the input order of candidate cuts is inconsistent. As stated in \cite{bixby1992implementing, maros2002computational, wang2023learning}, the order of constraints in linear programming (LP) affects the initial basis constructed. Adding different sequences of cuts to the LP relaxation also significantly impacts the efficiency of solving MILP problems. For example, we select a specific MILP problem from the CORLAT problem set~\cite{conrad2007connections,gomes2008connections} and use a trained sequence model for cut selection. We conduct six different experiments, with the first one using the default cut order and others randomly shuffling the input sequence (cuts) based on a random seed, shown in Table~\ref{motivation} in Appendix~\ref{MoreMotivation}. Then, according to these cut sequences, the MILP problem is solved respectively, and the experimental results are recorded in Figure~\ref{motivation_f}, including the solving time and primal-dual gap integral. From the experimental results, it is evident that different orders of cuts result in significant variations in both the solving time and primal-dual gap integral, thereby significantly impacting the stability of the solving process.

In summary, even if the candidate cuts and the MILP problem are fixed, the sequence models considering the positional information of the input sequence leads them to add different sequences of cuts to the LP relaxation based on the different orders of the input sequence. This significantly affects problem-solving. Additional experimental results regarding motivation can be found in Appendix~\ref{MoreMotivation}.

\section{Methodology}
\label{sec3}

\begin{figure}[h]
\centering
\includegraphics[width=0.85\textwidth]{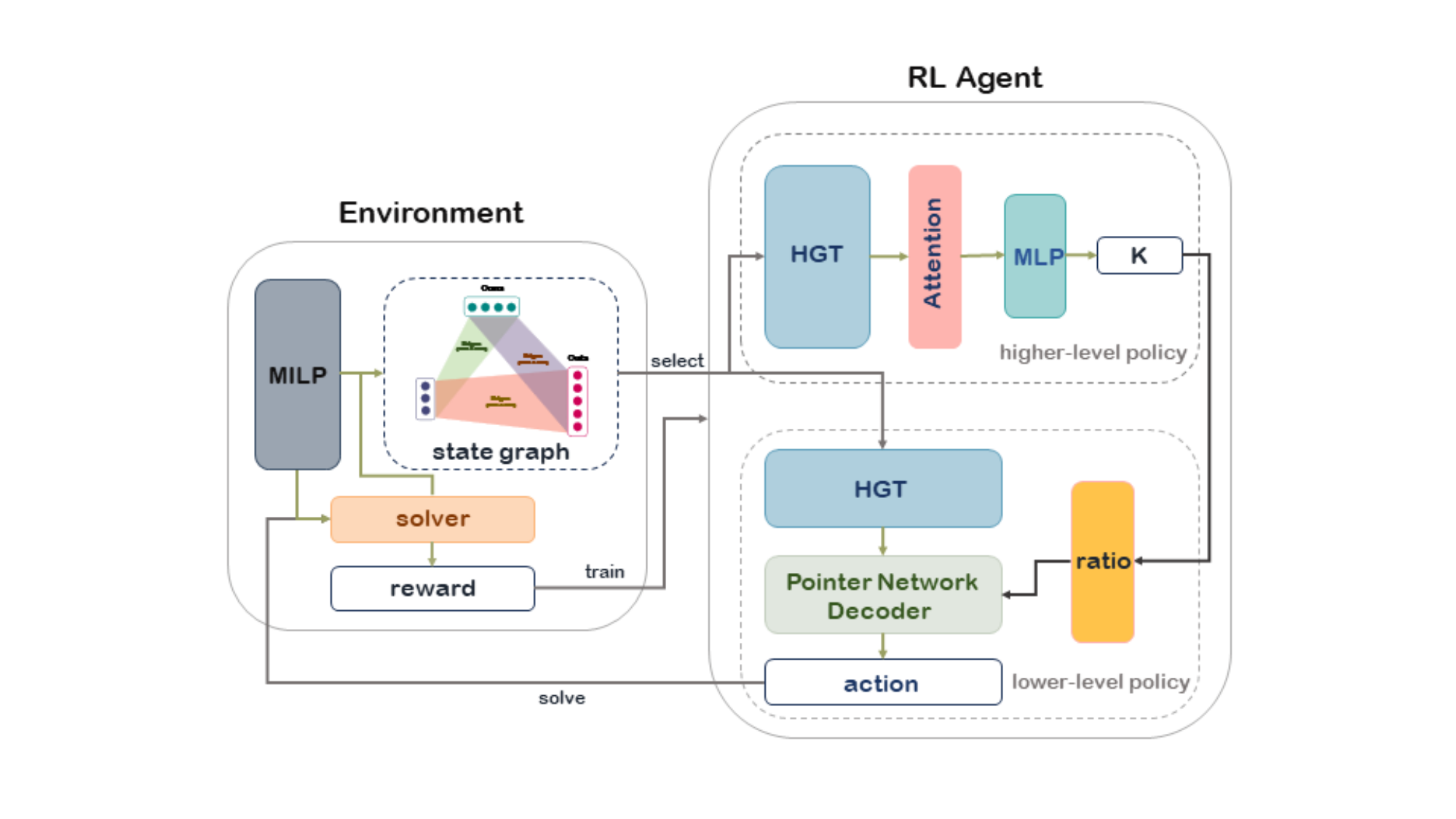}
\vspace{-20pt}
\caption{\footnotesize The framework of our model. SCIP acts as the environment, and the reinforcement learning model acts as the agent. SCIP acts as the environment, responsible for extracting MILP state features and solving the problem based on the cuts selected by the RL agent. The policy network, serving as the RL agent, encodes the state information into a tripartite graph. It is composed of two parts: the higher-level network, which outputs the proportion of selected cut planes, and the lower-level network, which generates the specific sequence of cuts.}
\label{structure_g}
\vspace{-10pt}
\end{figure}

In this section, we primarily introduce the specific structure of our model \textbf{HGTSM}. Inspired by the work~\cite{wang2023learning}, we formulate the cut selection problem as a Markov decision process, with adopting a similar hierarchical reinforcement learning framework for cut selection. Differently, we represent the MILP state as a tripartite graph and utilize a heterogeneous graph network to capture the underlying information of MILP. Subsequently, a modified sequence model is employed to generate the cut sequence. We emphasize that the new sequence model is more suitable for the task of cut selection. The framework of our model is shown in Figure~\ref{structure_g}. We use SCIP to obtain a tripartite graph representation of the problem state and use it as input for the policy network. The policy network consists of two networks: the higher-level network is responsible for calculating the number of selected cuts, and this output is then used as the input to the lower-level network to obtain the sequence of selected cuts.

\vspace{-10pt}

\subsection{Markov Decision Process}
\label{Markov}
As stated in \cite{tang2020reinforcement, wang2023learning}, the cut selection problem can be formalized as a Markov Decision Process (MDP). We define the MDP as a quadruple $(\mathcal{S},\mathcal{A},r,f)$, where $\mathcal{S}$ represents the state space, $\mathcal{A}$ represents the action space, $r:\mathcal{S}\times\mathcal{A}\rightarrow\mathbb{R}$ denotes the reward function, and $f$ represents the state transition function. At timestamp $t\neq0$, the environment is under state $s_t\in \mathcal{S}$, and the agent selects action $a_t\in \mathcal{A}$ based on the current state $s_t$. Subsequently, the environment returns a reward and transitions to the next state, repeating this process until an episode ends.

Our goal is to enable the agent to learn a policy $\pi$ that maximizes the cumulative reward. The policy $\pi:\mathcal{S}\rightarrow \mathcal{P}(\mathcal{A})$ provides a mapping from states to action distributions. Next, we formalize the cut selection problem as an MDP. We regard the solver as the environment and the deep reinforcement learning model as the agent, then specify the meanings of the state space $\mathcal{S}$, action space $\mathcal{A}$, reward function $r$, and state transition function $f$ as follows.

\textbf{State space $\mathcal{S}$.} The state space $\mathcal{S}$ represents the set of all possible states that the environment can be in. A state $s_t$ at timestamp $t$ is defined as $s_t = \{LP_t,C_t\}$. $LP_t$ represents the current LP relaxation problem, and $C_t=\{e_i^T\mathbf{x}\le d_i\}_{i=1}^{N_t}$ $(e_i\in \mathbb{R}^n,d_i\in \mathbb{R})$ represents the current set of candidate cuts, where $N_t$ is the number of cuts. We elaborate on the encoding of $LP_t$ and $C_t$ in Section~\ref{StateEncoding}. 

\textbf{Action space $\mathcal{A}$.} The action space $\mathcal{A}$ represents the set of all possible actions that the agent can take under the state $s_t \in S$. Specifically for the cut selection problem, assuming the current solver is in state $s_t$ with the set of candidate cuts $C_t$, the action space $\mathcal{A}$ is defined as the set of all ordered subsets of $C_t$, and a specific action $a_t$ is then an ordered subset of $C_t$.

\textbf{Reward function $r$.} The reward function $r(s_t, a_t)$ represents the feedback from the environment regarding the impact of action $a_t$ taken by the agent under state $s_t$. To assess the effect of adding cuts on problem-solving, we define the reward as the solving time or the primal-dual gap integral and collect reward data after the problem-solving process ends.

\textbf{State transition function $f$.} The state transition function describes how the environment transitions from state $s_t$ to the next state $s_{t+1}$ after the agent takes action $a_t$ under state $s_t$. Specifically for this problem, $s_{t+1}$ represents the new $LP_{t+1}$ relaxation and the new candidate set of cuts $C_{t+1}$ generated by the solver under the current $LP_t$ relaxation.

\subsection{State Encoding}
\label{StateEncoding}
\begin{figure}[H]
\centering
\subfigure[ ]{
\includegraphics[width=0.45\textwidth]{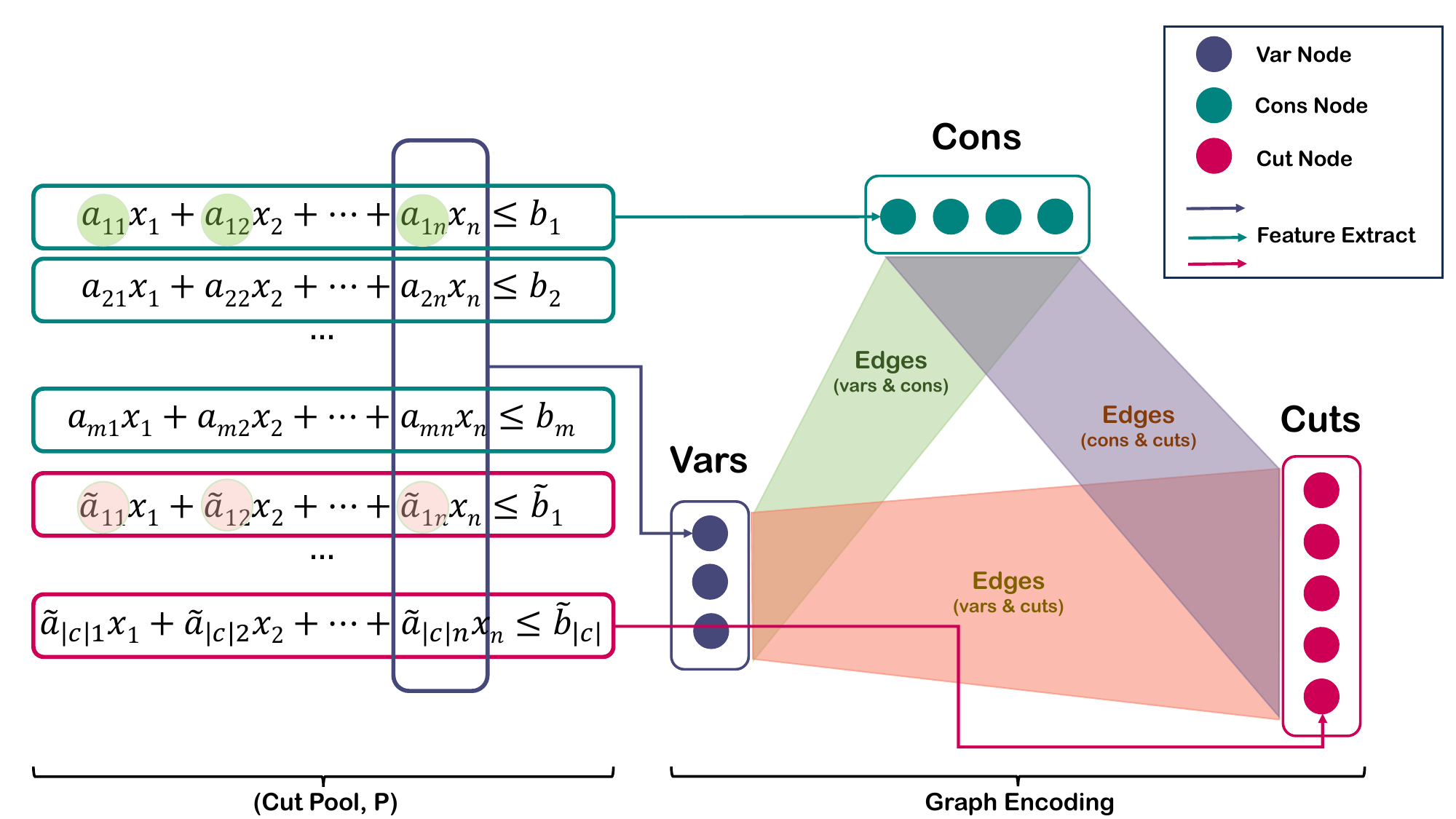}}
\subfigure[ ]{
\includegraphics[width=0.45\textwidth]{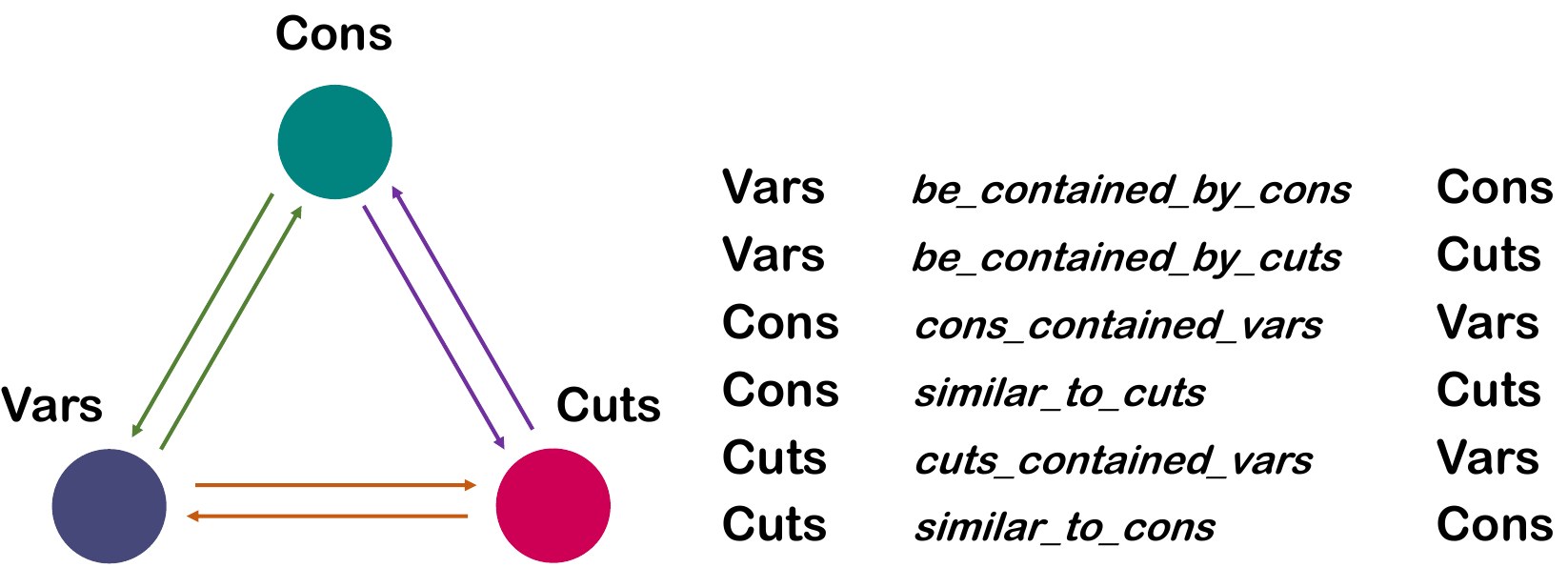}}
\caption{(a) Encode a state $s$ into a heterogeneous tripartite graph. Constructing edges between constraints (or cuts) and variables with non-zero coefficients and between constraints with non-zero parallelism coefficients and cuts. (b) The schema and meta relations of state graph.}
\label{Fig.main}
\vspace{-10pt}
\end{figure}

As mentioned in Section~\ref{Markov}, we represent the state $s$ of the solver at timestamp 
$t$ using the LP relaxation $LP_t$ and the candidate cuts $C_t$. Considering that the LP relaxation consists of the objective function $\mathbf{c}^T\mathbf{x}$ and constraints $\mathbf{A}\mathbf{x}\le \mathbf{b}$, where each constraint $a_ix\leq b_i$ is associated with different variables, we manually construct features for constraints ($Cons$) and variables ($Vars$) to characterize the LP relaxation. Additionally, we construct features for cuts as the representation of $Cuts$. To better characterize the solver's state and establish the inherent connections among $Cons$, $Vars$, and $Cuts$, we adopt a similar representation in ~\cite{gasse2019exact, paulus2022learning} and represent the state $s$ as a heterogeneous tripartite graph. As shown in Figure~\ref{Fig.main}(a),  three types of nodes represent $Cons$, $Vars$, and $Cuts$, respectively. The edges between nodes  represent the relationships among $Cons$, $Vars$, and $Cuts$. 

Different from the aforementioned graph representation, we no longer directly use the coefficients of variables in constraints or candidate cuts to establish relationships between nodes. Instead, we obtain attention scores between different nodes through attention mechanisms to align with the heterogeneous graph network.
Furthermore, different from the tripartite graph as an undirected graph~\cite{gasse2019exact,labassi2022learning,paulus2022learning}, we represent the relationships between nodes as bidirectional directed graphs, denoted as $G=\{\mathcal{V},\mathcal{E},\mathcal{A},\mathcal{R}\}$, as shown in Figure~\ref{Fig.main}(b), where  $\mathcal{V}$ is the node set, $\mathcal{E}$ is the edge set, $\mathcal{A}$ represents the node types, and $\mathcal{R}$ represents the edge types. A directed edge $e=(s,t)$ means the edge from source node $s$ to target node $t$. It is remarkable that the relationships from $Vars$ to $Cons$ are not the same as the relations from $Cons$ to $Vars$. 
We believe that considering both opposite relationships are helpful to extract state information comprehensively. Similar designs are adopted for relationships between other types of nodes.

The next step is to construct feature vectors. 
For $Vars$ and $Cons$, we use the feature design proposed by Paulus et al.~\cite{paulus2022learning}. As to $Cuts$, we adopt the design proposed by Wang et al.~\cite{wang2023learning}. In summary, the state we construct is represented by features $Vars\in R^{n\times17}$, $Cons\in R^{m\times16}$, and $Cuts\in R^{c\times13}$, where $n$ is the number of $Vars$, $m$ is the number of original $Cons$, and $c$ is the number of $Cuts$. Table~\ref{feature} in Appendix \ref{FeatureDetails} provides a detailed description of the features.

\subsection{Heterogeneous Graph Neural Network}
To leverage the tripartite graph constructed in Section ~\ref{StateEncoding}, different from the works of ~\cite{gasse2019exact, labassi2022learning, paulus2022learning}, we utilize heterogeneous graph neural networks to extract different node information and inter-node relationships. 
Because graph convolutional networks (GCN) used in the aforementioned works do not fully exploit the characteristics of heterogeneous graphs. For example, there are two message passing steps, $Vars \rightarrow Cons$ and $Cons \rightarrow Cuts$, where $Vars$ and $Cons$ serve as source nodes, and $Cons$ and $Cuts$ act as target nodes. However, in GCN, all source nodes undergo the same network processing, whatever these source nodes are in different types. Similarly, all target nodes and edges also undergo the same network processing. Obviously, GCNs treat different nodes and edges with the same processing but neglects the properties of heterogeneous graphs, which may cause insufficient extraction of state information in solvers.

To extract more information from heterogeneous graphs, we are inspired by \cite{hu2020heterogeneous} to utilize a Heterogeneous Graph Transformer (HGT), which can construct network structures specific to different types of nodes or edges, and preserve the heterogeneity of the graph. By computing the attention scores of source nodes to target nodes, corresponding weights are obtained for message propagation and aggregation and used to update the feature vectors of each node, as shown in Eq.~(2). 
\begin{equation}
H^{(l)}[t]=linear_{\tau(t)}(\sigma(\widetilde{H}^{(l)}[t]))+H^{(l-1)}[t], 
\end{equation}
where $\widetilde{H}^{(l)}[t]=\mathop{\oplus}\limits_{\forall s\in N(t)}(Attention_{HGT}(s,e,t)\cdot Message_{HGT}(s,e,t)).$

Here, $H^{(l)}[t]$ represents the output of the network at the $l\mbox{-}th$ layer and serves as the new feature of the target node. $\sigma$ denotes the ReLU activation function, and $linear_{\tau(t)}$ is a linear layer specific to the node type. $\widetilde{H}^{(l)}[t]$ represents the updated vector aggregated with neighborhood information, $N(t)$ denotes the set of source nodes for node $t$, $e$ represents the edge from source node $s$ to target node $t$, $Attention$ computes the importance of source nodes, and $Message$ extracts message from source nodes. Appendix \ref{HGT} provides more details about the Heterogeneous Graph Transformer. After message propagation and aggregation, we discard the features of $Vars$ and $Cons$, retaining only $Cuts$ as the main features since it includes all information of the LP relaxation. 

\subsection{Policy Network}
As Wang et al.~\cite{wang2023learning} pointed out, the core of the cut selection problem lies in the number of selected cuts and the order in which these cuts should be selected. Therefore, based on ~\cite{wang2023learning}, we design a hierarchical reinforcement learning approach to address the cut selection problem. In this approach, the higher-level policy network $\pi^h:\mathcal{S}\rightarrow\mathcal{P}([0,1])$ is utilized to compute the ratio of selected cuts from the candidate set, modeled as a tanh-Gaussian distribution. The mean and variance of this distribution are inferred by neural networks, and the results are then mapped to the interval $[0,1]$ through linear transformation. The higher-level policy network is parameterized as follows:
\[
\pi_{\theta_h}^h(\cdot|s) = 0.5\times tanh(K) + 0.5, K\sim \mathcal{N}(\mu_{\theta_h}(s),\sigma_{\theta_h}(s)). 
\]
Here, $\mu_{\theta_h}(s)$ and $\sigma_{\theta_h}(s)$ represent the mean and standard deviation of the neural network output in state $s$, respectively. $K$ is the result of random sampling from a Gaussian distribution. The lower-level policy network $\pi^l:\mathcal{S}\times[0,1]\rightarrow\mathcal{P}(\mathcal{A})$, based on the output of the higher-level network, executes specific cut selection actions. Specifically, suppose the candidate set of cuts at state $s_t$ is $C_t$, with a total of $|C_t|$ cuts. If the output of the higher-level policy network is $ratio$, then the lower-level policy network needs to select a sequence of cuts from the candidate set of cuts, with the number of selected cuts being $|C_t|\times ratio$. This implies that the cut selection process can be modeled as a sequence-to-sequence problem where the elements of the output sequence are entirely drawn from the input sequence. This characteristic aligns well with pointer networks. Hence, the lower-level policy network primarily utilizes pointer networks to execute the specific cut selection process.

In terms of state representation, previous reinforcement learning-based approaches often choose to use the features of candidate cuts as inputs to the model. While models like pointer networks with attention mechanisms can incorporate information from other cuts in the sequence, we believe that considering only adjacent cuts is insufficient; the model needs to thoroughly consider the states of the solver. Therefore, we extract features for $Cuts$, $Vars$, and $Cons$ separately to construct a tripartite graph as the input to the model. After aggregating information through HGT, we obtain the features of $Cuts$. Subsequently, we calculate $ratio$ and select sequences accordingly. Moreover, another advantage of graph networks is their ability to handle MILP problems of different scales.
On the other hand, we argue that the input order of candidate cuts should not influence the cut selection. Hence, neither the calculation of the ratio by the higher-level policy network nor the selection of sequences by the lower-level network should consider the positional information of the input sequence of cuts because the cut order does not reflect the geometry of the feasible set~\cite{tang2020reinforcement}. In other words, we aim for the network to possess permutation invariance.

Based on this idea, the policy network architecture we construct is illustrated in the right of Figure~\ref{structure_g}. The processing of the higher-level network is as follows: the input features encoded as a tripartite graph are first processed by HGT. The processed features of $Vars$ and $Cons$ are then discarded, leaving behind the aggregated features of $Cuts$. Subsequently, attention mechanism is employed to comprehensively consider the features of all candidate $Cuts$. Finally, a multilayer perceptron (MLP) is employed to obtain the mean and variance of a Gaussian distribution, from which the ratio of selected cuts is determined. The lower-level network modifies the traditional pointer network, by using HGT as the encoder part of the pointer network. Based on the ratio output by the higher-level network, the processed features of candidate cuts are passed through the decoder part of the pointer network to obtain the selected sequence of cuts.

To the best of our knowledge, we are the first to combine heterogeneous graph neural networks with deep reinforcement learning algorithms to address the cut selection problem. Unlike previous learning-based approaches, the use of graph neural networks enables the model to not only consider the information of $Cuts$ but also fully consider the relationship between $Cuts$ and LP relaxations, thus aiding the model in capturing potential information more effectively. Additionally, deep reinforcement learning algorithms allow the model to be trained without pre-labeled information, facilitating the exploration of more efficient cut selection strategies.

\subsection{Training}
For policy-based reinforcement learning algorithms, the goal is to maximize the objective function:
\[
J(\theta) = \mathbb{E}_{s\sim \mu, a_t\sim \pi_{\theta}(\cdot | s)}[Q_\pi (s,a_t)],
\]
where $\mu$ denotes the initial state distribution, $\theta = [\theta_h, \theta_l]$ refers to the parameters of the higher-level network and the lower-level network, respectively. $a_t$ indicates the action selected by the policy network $\pi$ under the state $s$, and $Q_\pi(s,a_t)$ represents the action-value function, which is used to measure the value of action $a_t$ under state $s$. Since this paper employs hierarchical reinforcement learning~\cite{pateria2021hierarchical} to find the parameters $\theta$ that maximize the objective function, 
the policy is calculated with $\pi_{\theta}(\cdot | s) = E_{p\sim \pi_{\theta_h}^h(\cdot | s)}[\pi_{\theta_l}^l(a_t | s, ratio)]$. Based on the policy gradient theorem:
\[
\nabla_{\theta_h}J([\theta_h,\theta_l]) = \mathbb{E}_{s\sim \mu, ratio\sim \pi_{\theta_h}^h(\cdot | s)}[\nabla_{\theta_h}\log(\pi_{\theta_h}^h(ratio|s))\mathbb{E}_{a_t\sim \pi_{\theta_l}^l(\cdot| s,ratio)}[Q_\pi (s,a_t)]],
\]
\[
\nabla_{\theta_l}J([\theta_h,\theta_l]) = \mathbb{E}_{s\sim \mu, ratio\sim \pi_{\theta_h}^h(\cdot | s), a_t\sim \pi_{\theta_l}^l(\cdot| s,ratio)}[\nabla_{\theta_l}\log\pi_{\theta_l}^l(a_t|s,ratio)Q_\pi (s,a_t)],
\]
we train the policy network using the Actor-Critic~\cite{mnih2016asynchronous} algorithm with a baseline. 

\section{Experiments}
\label{sec4}
To better evaluate our model, we mainly conduct experiments from three aspects, listed as follows:
\begin{itemize}
\item Experiment 1: \textbf{Performance evaluation}. We will evaluate our model on six challenging sets of MILP problems and compare it with SOTA learning-based methods. 
\item Experiment 2: \textbf{Stability evaluation}. We will conduct comparative experiments with the SOTA sequence-based model to evaluate the stability of our proposed model in problem-solving.
\item Experiment 3: \textbf{Generalization evaluation}. We will test the generalization of our model. Be specific, the model is  trained on one dataset, and then evaluate it on other datasets with  different MILP problem types.
\end{itemize}

\subsection{Benchmarks}
Our benchmarks are based on the MILP problem sets used in \cite{wang2023learning}, from which we select six challenging NP-hard problem sets. Following the classification method of \cite{wang2023learning}, we categorize the problem sets into two types based on the solving time of SCIP: medium and hard. Note that we omit the easy dataset in \cite{wang2023learning} because its MILP problems can be easily solved in one minute. The medium difficulty problem sets include MIK~\cite{atamturk2003facets} and Corlat~\cite{conrad2007connections,gomes2008connections}, while the hard difficulty problem sets include Load Balancing, Anonymous~\cite{ecoleai2021}, MIPLIB mixed neos and MIPLIB mixed supportcase. Detailed information about the MILP problem datasets can be found in Appendix~\ref{DatasetsDetails}.

\subsection{Baselines}
The baselines in our comparative experiments include five heuristic cut selection methods: \textit{No Cuts, Random, Normalized Violation (NV), Efficacy (Eff)}, and \textit{the Default method} in SCIP 8.0.0. We also compare with the SOTA learning-based selection method HEM \cite{wang2023learning}. More descriptions of these baselines can be found in Appendix \ref{Baselines}.

\subsection{Evaluation metrics}
To evaluate the problem-solving capability of our model, we select two widely used metrics: average solving time and average primal-dual gap integral. These metrics serve as evaluation indicators for the performance, assessing both the efficiency and quality of the solved problems. Additionally, to measure the stability of our model in solving problems with different types, we propose a new metric $Stability$ as an evaluation criterion, which is calculated as follows,  
\[
Stability = \{ \frac{\sum_{i=1}^{N}(\max(IR(p),ID(p))/\min(IR(p),ID(p)))}{N} - 1, p \in D\}
\]
Here, $IR(p)$ represents the set of primal-dual integrals for each problem in datasets after shuffling the order of input cuts with the specified random seed. $ID(p)$ represents the set of primal-dual integrals after inputting the model with the default sequence of cuts generated by SCIP. $D$ represents the specified MILP datasets. $N$ represents the number of problems in the dataset. $Stability$ calculates the offset ratio of the integral under different random seeds while avoiding inaccurate stability evaluation caused by mutual interference between multiple problems, and effectively evaluating the stability of the model for solving MILP problems. Additionally, we introduce the metric $ Improvement$, which represents the enhancement of our model's stability compared to HEM~\cite{wang2023learning}. $Improvement$ is calculated as follows. 
\[
    Improvement =  \frac{S(HEM) - S(HGTSM)}{S(HEM)}.
\]
Here, $S(\cdot)$ represents the $Stability$ of the model on a given dataset.

\subsection{Experiments and Results}
\paragraph{Experiment 1: Performance evaluation} We evaluate the baselines and our model on six challenging MILP datasets, with the results shown in Table~\ref{result_1_1}. The results demonstrate that our model outperforms the baselines significantly on both the medium and hard datasets. Particularly, our model achieves a significant advantage over the heuristic baselines on the primal-dual gap integral. Moreover, we also substantially reduce the solving time on the medium dataset. Even compared to the SOTA method HEM, our model still maintains a clear advantage. Remark that the SOTA method we  fine-tuned not only outperforms the heuristic methods significantly but also performs even better than the original paper \cite{wang2023learning}. This further enhances the persuasiveness of our model.

\begin{table}[h]
\caption{The statistical data on six datasets include the mean (standard deviation) of solving time and primal-dual gap integral. }
\label{result_1_1}
\centering
\scalebox{0.73}{
\begin{tabular}{lllllll}
\toprule
\multicolumn{3}{c}{Medium: MIK} & \multicolumn{2}{c}{Medium: Corlat} & \multicolumn{2}{c}{Hard: Load Balancing} \\
\cmidrule(r){1-3} \cmidrule(r){4-5} \cmidrule(r){6-7}
Method & Time(s) $\downarrow$ & PD Integral $\downarrow$ & Time(s) $\downarrow$ & PD Integral $\downarrow$ & Time(s) $\downarrow$ & PD Integral $\downarrow$\\
\cmidrule(r){1-3} \cmidrule(r){4-5} \cmidrule(r){6-7}
Nocuts & 300.01(0.002) & 2150.49(1041.24) & 93.95(125.14) & 2494.07(6160.41) & 300.03(0.07) & 14785.51(948.62) \\
Default & 161.56(127.59) & 769.49(899.79) & 80.05(122.17) & 2545.65(6073.38) & 300.02(0.06) & 11934.16(2691.11) \\
Random & 266.83(66.40) & 1829.29(941.69) & 77.52(118.78) & 2141.34(5866.99) & 300.07(0.20) & 13687.07(1215.45) \\
Nvscore & 289.99(30.04) & 2355.49(636.03) & 83.65(126.16) & 3072.78(7098.79) & 300.02(0.04) & 13971.28(962.97) \\
Effscore & 291.32(18.289) & 2237.09(662.46) & 93.00(124.91) & 2610.05(6017.46) & 300.04(0.09) & 13957.91(966.96) \\
\cmidrule(r){1-3} \cmidrule(r){4-5} \cmidrule(r){6-7}
HEM & 162.02(122.82) & 786.53(774.79) & 56.70(109.67) &  933.86(2141.39) & 300.02(0.02) & 9612.05(1030.69) \\
HGTSM & \textbf{156.60(129.93)} & \textbf{669.47(657.30)} & \textbf{56.23(107.51)} & \textbf{789.31(1706.95)} & \textbf{300.01(0.02)} &  \textbf{9530.80(1051.30)} \\
\midrule
\multicolumn{3}{c}{Hard: Anonymous} & \multicolumn{2}{c}{Hard: MIPLIB mixed neos} & \multicolumn{2}{c}{Hard: MIPLIB mixed supportcase} \\
\cmidrule(r){1-3} \cmidrule(r){4-5} \cmidrule(r){6-7}
Method & Time(s) $\downarrow$ & PD Integral $\downarrow$ & Time(s) $\downarrow$ & PD Integral $\downarrow$ & Time(s) $\downarrow$ & PD Integral $\downarrow$\\
\cmidrule(r){1-3} \cmidrule(r){4-5} \cmidrule(r){6-7}
Nocuts & 244.84(97.09) & 17976.51(9706.82) & 248.78(93.91) & 13956.90(12927.56) & \textbf{156.05(132.41)} & 9692.81(10764.72) \\
Default & 239.70(99.69) & 16399.05(9299.57) & 251.03(84.82) & 14133.80(12636.66) & 161.64(135.20) & 9689.05(11339.05) \\
Random & 240.05(100.40) & 16825.56(9236.02) & 247.23(91.40) & 13990.47(12653.49) & 164.57(132.04) & 9786.02(10881.72) \\
Nvscore & 239.70(99.69) & 16399.05(9299.57) & 254.70(78.46) & 14538.32(12378.39) & 167.08(129.24) & 9095.64(10138.13) \\
Effscore & 241.02(100.61) & 17423.81(9932.39) & 247.72(94.02) & 13824.76(13035.09) & 169.41(128.01) & 8761.83(8985.01) \\
\cmidrule(r){1-3} \cmidrule(r){4-5} \cmidrule(r){6-7}
HEM & \textbf{238.10(107.61)} & 15217.53(9204.87) & 247.20(91.44) & 8419.75(12479.65) & 163.30(136.95) & 7188.17(6875.03) \\
HGTSM & 241.07(104.20) & \textbf{13032.58(9386.95)} & \textbf{246.44(89.29)} & \textbf{8326.77(12514.00)} & 164.46(136.77) & \textbf{5644.90(5893.50)} \\
\bottomrule
\end{tabular}
}
\end{table}

\begin{figure}[H]
\centering
\subfigure[ ]{
\includegraphics[width=0.43\textwidth]{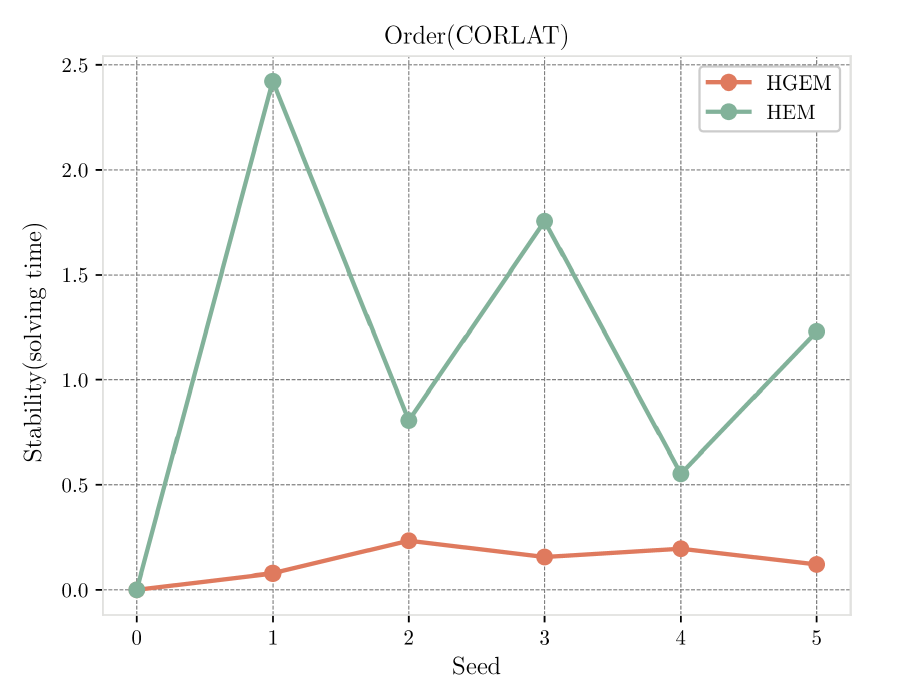}}
\subfigure[ ]{
\includegraphics[width=0.43\textwidth]{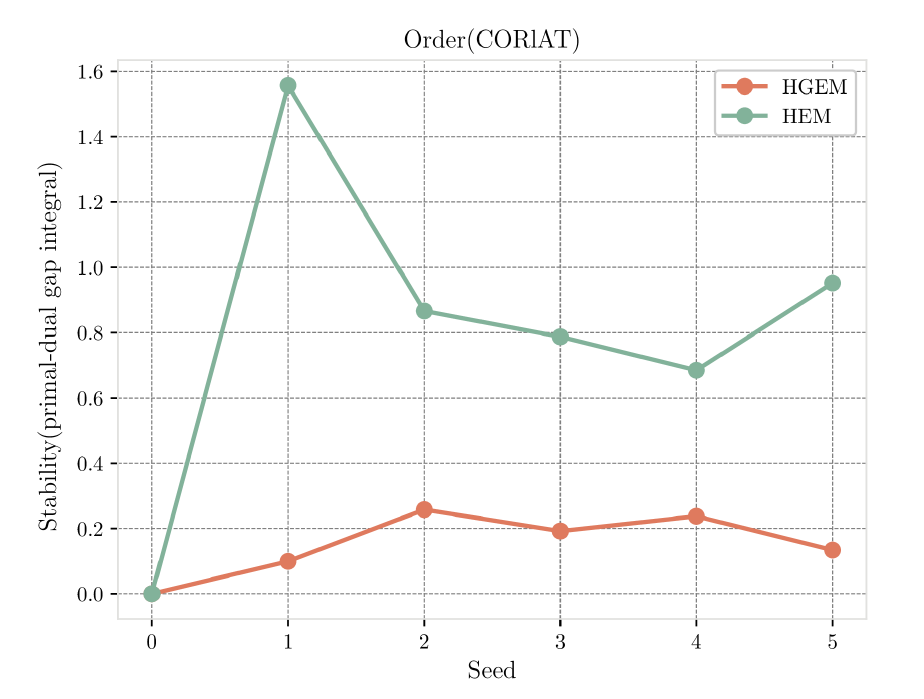}}
\caption{(a) $Stability$ scores based on solving time are calculated using different random seeds. (b) $Stability$ scores based on the primal-dual gap integral are calculated using different random seeds. In both figures, the horizontal ordinate $0$ represents the use of the default input sequence from SCIP without modifications.}
\label{sta_result}
\end{figure}

\vspace{-10pt}
\paragraph{Experiment 2: Stability evaluation}
\label{ExpSta}
As mentioned above, we notice that sequence models capture the positional information of the input sequence, which is detrimental to the solving stability. Therefore, we improve the sequence model by utilizing the permutation invariance of graph networks to ignore the positional information for enhancing the solving stability. We randomly shuffle the cut sequences with different random seeds and calculate the $Stability$ scores (the lower, the better) in terms of solving time and primal-dual gap integral. We report the trends in solving time and the stability of the primal-dual gap integral on the CORLAT dataset under different random seeds, as shown in Figure~\ref{sta_result}. 
Compared to the SOTA method HEM~\cite{wang2023learning}, our model shows significant advantages in both the stability of solving time and the stability of solution quality. 
Additionally, as shown in Table~\ref{result_2_1} and Table~\ref{result_2_2}, we obtain the $Stability$ score for each dataset by averaging the $Stability$ scores across different random seeds. We then present the $Improvement$ of our model's stability compared to HEM in terms of solving time and PD integral.
Since the solving times for all problems in the Load balancing test set exceed $300$ seconds, and the MIPLIB mixed neos test set had too few problems, resulting in significant serendipity, we exclude these two datasets and report the experimental results on other datasets in Table~\ref{result_2_1} and Table~\ref{result_2_2}.
The experimental results show that compared to HEM, our model exhibits significant stability improvements across multiple datasets. Notably, on the CORLAT dataset, the improvements in both solving time and the primal-dual gap integral exceed $80\%$.
These demonstrate the importance of ignoring the positional information of the cut sequence. More results regarding the stability evaluation can be found in the Appendix \ref{MoreStaResult}.

\paragraph{Experiment 3: Generalization evaluation.}
We hope the model based on graph networks can capture more general underlying information of MILP problems. Therefore, we evaluate the generalization of our model across different categories of MILP problems. We train our model on one dataset and then evaluate the model on other MILP datasets. The experimental results regarding the primal-dual gap integral are shown in Table~\ref{result_Gen}. According to the results, the models trained for specific problems exhibit the best performance for those problems (the diagonal line in Table~\ref{result_Gen}). However, the models also demonstrate acceptable performance on other problems. Although the models are not trained for certain problems, surprisingly, their performance generally exceeds that of heuristic methods and even outperforms learning-based baselines in some cases. The experimental results clearly demonstrate the potentiality and prospect of finding more general and efficient methods for solving MILP problems. More generalization results of our model are provided in Appendix \ref{MoreGenResult}.

\begin{table}
\caption{The statistical comparison on five datasets includes the mean (standard deviation) of the primal-dual gap integral. Note that Sta. means stability, and PD Imp. means the improvement(\%) of PD integral.}
\label{result_2_1}
\centering
\scalebox{0.7}{
\begin{tabular}{lllllllllll}
\toprule
~ & \multicolumn{2}{c}{CORLAT} & \multicolumn{2}{c}{MIK} & \multicolumn{2}{c}{Anonymous} & \multicolumn{2}{c}{MIPLIB mixed supporecase} & \multicolumn{2}{c}{Load Balancing}\\
\cmidrule(r){1-3} \cmidrule(r){4-5} \cmidrule(r){6-7} \cmidrule(r){8-9} \cmidrule(r){10-11}
Method  & Sta. $\downarrow$ & PD Imp. $\uparrow$ & Sta. $\downarrow$ & PD Imp. $\uparrow$ & Sta. $\downarrow$ & PD Imp. $\uparrow$ & Sta. $\downarrow$ & PD Imp. $\uparrow$ & Sta. $\downarrow$ & PD Imp. $\uparrow$ \\
\cmidrule(r){1-3} \cmidrule(r){4-5} \cmidrule(r){6-7} \cmidrule(r){8-9} \cmidrule(r){10-11}
HEM & 0.96 & NA & 0.41 & NA & 0.44 & NA & 0.44 & NA & 0.02 & NA\\
HGTSM & \textbf{0.18} & \textbf{81.02} & \textbf{0.14} & \textbf{65.46} & \textbf{0.25} & \textbf{43.61} & \textbf{0.18} & \textbf{57.94} & \textbf{0.004} & \textbf{81.03} \\
\bottomrule
\end{tabular}
}
\end{table}

\begin{table}
\caption{The generalization results of our model are based on the metric of primal-dual gap integral. For example, $1193.89(3034.38)$ in the cell $[2,1]$ of datagrid means that the model is trained on the dataset MIK and evaluated on the dataset CORLAT. $1193.89$ is the average of primal-dual integrals, and $3034.38$ is the standard deviation of primal-dual integrals. Note that Neos means MIPLIB mixed neos, and Supportcase means MIPLIB mixed supportcase. }
\label{result_Gen}
\centering
\scalebox{0.66}{
\begin{tabular}{lllllll}
\toprule
~ & \multicolumn{6}{c}{Generalization results of Primal-dual gap integral} \\
\cmidrule(r){2-7}
~ & CORLAT & MIK & Anonymous & Neos & Supportcase & Load Balancing \\
\midrule
CORLAT & \textbf{789.31(1706.95)} & 1357.02(855.48) & 13888.85(9321.71) & 11965.77(11719.81) & 8194.37(9296.81) & 11163.01(962.70) \\
MIK & 1193.89(3034.38) & \textbf{669.47(657.30)} & 14674.28(9065.02) & 11971.79(11332.99) & 9564.82(10878.71) & 9790.31(1036.44)\\
Anonymous & 1574.03(4449.62) & 1259.07(830.05) & \textbf{13032.58(9386.95)} & 13681.84(12524.28) & 7917.94(9502.84) & 12175.21(736.33)\\
Neos & 3274.90(7284.15) & 1166.83(857.64) & 13832.84(8755.20) & \textbf{8326.77(12514.00)} & 10740.65(11788.25) & 11688.87(1006.33)\\
Supportcase & 2078.57( 5752.44) & 1579.84(880.53) & 16238.75(9304.63) & 13865.66(13088.69) & \textbf{5644.90(5893.50)} & 12139.26(792.47) \\
Load Balancing & 1820.51(5870.83) & 1002.38( 828.96) & 15366.35(8413.04) & 13699.30(12644.24) & 9648.57(11816.80) & \textbf{9530.80(1051.30)}\\
\bottomrule
\end{tabular}
}
\end{table}

\section{Conclusion}
\label{sec5}
Cut selection is crucial for solving MILP problems. This paper proposes a novel reinforcement learning model to learn the cut selection strategy. We encode the state as a tripartite graph and utilize the HGT to better capture the underlying information of MILP. We emphasize that the order of cuts inputs should not affect the solving efficiency. Therefore, we use a modified pointer network as a sequence model to generate the cut sequence. Compared to heuristic methods and SOTA learning-based method, the proposed model further improves the performance of MILP solving and significantly improves the stability of problem solving. Additionally, we discuss the generality of our model for different types of problems. For future work, we hope that our model can be effective under more complex conditions and further get improvements in solving time. Exploring more general MILP-solving strategies would be an interesting and meaningful topic.

\bibliographystyle{plain}

\bibliography{references}

\begin{thebibliography}{10}

\bibitem{atamturk2003facets}
Alper Atamt{\"u}rk.
\newblock On the facets of the mixed--integer knapsack polyhedron.
\newblock {\em Mathematical Programming}, 98(1):145--175, 2003.

\bibitem{bengio2021machine}
Yoshua Bengio, Andrea Lodi, and Antoine Prouvost.
\newblock Machine learning for combinatorial optimization: a methodological tour d’horizon.
\newblock {\em European Journal of Operational Research}, 290(2):405--421, 2021.

\bibitem{bestuzheva2021scip}
Ksenia Bestuzheva, Mathieu Besan{\c{c}}on, Wei-Kun Chen, Antonia Chmiela, Tim Donkiewicz, Jasper van Doornmalen, Leon Eifler, Oliver Gaul, Gerald Gamrath, Ambros Gleixner, et~al.
\newblock The scip optimization suite 8.0.
\newblock {\em arXiv preprint arXiv:2112.08872}, 2021.

\bibitem{bixby1992implementing}
Robert~E Bixby.
\newblock Implementing the simplex method: The initial basis.
\newblock {\em ORSA Journal on Computing}, 4(3):267--284, 1992.

\bibitem{conrad2007connections}
Jon Conrad, Carla~P Gomes, Willem-Jan Van~Hoeve, Ashish Sabharwal, and Jordan Suter.
\newblock Connections in networks: Hardness of feasibility versus optimality.
\newblock In {\em Integration of AI and OR Techniques in Constraint Programming for Combinatorial Optimization Problems: 4th International Conference, CPAIOR 2007. Proceedings 4}, pages 16--28. Springer, 2007.

\bibitem{dey2018theoretical}
Santanu~S Dey and Marco Molinaro.
\newblock Theoretical challenges towards cutting-plane selection.
\newblock {\em Mathematical Programming}, 170:237--266, 2018.

\bibitem{ecoleai2021}
Maxime Gasse, Simon Bowly, Quentin Cappart, Jonas Charfreitag, Laurent Charlin, Didier Ch{\'{e}}telat, Antonia Chmiela, Justin Dumouchelle, Ambros~M. Gleixner, Aleksandr~M. Kazachkov, Elias~B. Khalil, Pawel Lichocki, Andrea Lodi, Miles Lubin, Chris~J. Maddison, Christopher Morris, Dimitri~J. Papageorgiou, Augustin Parjadis, Sebastian Pokutta, Antoine Prouvost, Lara Scavuzzo, Giulia Zarpellon, Linxin Yang, Sha Lai, Akang Wang, Xiaodong Luo, Xiang Zhou, Haohan Huang, Sheng~Cheng Shao, Yuanming Zhu, Dong Zhang, Tao Quan, Zixuan Cao, Yang Xu, Zhewei Huang, Shuchang Zhou, Binbin Chen, Minggui He, Hao Hao, Zhiyu Zhang, Zhiwu An, and Kun Mao.
\newblock The machine learning for combinatorial optimization competition {(ML4CO):} results and insights.
\newblock In Douwe Kiela, Marco Ciccone, and Barbara Caputo, editors, {\em NeurIPS 2021 Competitions and Demonstrations Track}, volume 176 of {\em Proceedings of Machine Learning Research}, pages 220--231. {PMLR}, 2021.

\bibitem{gasse2019exact}
Maxime Gasse, Didier Ch{\'e}telat, Nicola Ferroni, Laurent Charlin, and Andrea Lodi.
\newblock Exact combinatorial optimization with graph convolutional neural networks.
\newblock {\em Advances in Neural Information Processing Systems}, pages 15554--15566, 2019.

\bibitem{gleixner2021miplib}
Ambros Gleixner, Gregor Hendel, Gerald Gamrath, Tobias Achterberg, Michael Bastubbe, Timo Berthold, Philipp Christophel, Kati Jarck, Thorsten Koch, Jeff Linderoth, et~al.
\newblock Miplib 2017: data-driven compilation of the 6th mixed-integer programming library.
\newblock {\em Mathematical Programming Computation}, 13(3):443--490, 2021.

\bibitem{gomes2008connections}
Carla~P Gomes, Willem-Jan Van~Hoeve, and Ashish Sabharwal.
\newblock Connections in networks: A hybrid approach.
\newblock In {\em Integration of AI and OR Techniques in Constraint Programming for Combinatorial Optimization Problems: 5th International Conference, CPAIOR 2008 Proceedings 5}, pages 303--307. Springer, 2008.

\bibitem{gomory1960algorithm}
Ralph Gomory.
\newblock {\em An Algorithm for the Mixed Integer Problem: Notes Linear Programming and Extensions}.
\newblock Rand, s1960.

\bibitem{gupta2020hybrid}
Prateek Gupta, Maxime Gasse, Elias Khalil, Pawan Mudigonda, Andrea Lodi, and Yoshua Bengio.
\newblock Hybrid models for learning to branch.
\newblock {\em Advances in Neural Information Processing Systems}, pages 18087--18097, 2020.

\bibitem{he2014learning}
He~He, Hal Daume~III, and Jason~M Eisner.
\newblock Learning to search in branch and bound algorithms.
\newblock {\em Advances in Neural Information Processing Systems}, pages 3293--3301, 2014.

\bibitem{hu2020heterogeneous}
Ziniu Hu, Yuxiao Dong, Kuansan Wang, and Yizhou Sun.
\newblock Heterogeneous graph transformer.
\newblock In {\em Proceedings of the Web Conference 2020}, pages 2704--2710, 2020.

\bibitem{huang2022learning}
Zeren Huang, Kerong Wang, Furui Liu, Hui-Ling Zhen, Weinan Zhang, Mingxuan Yuan, Jianye Hao, Yong Yu, and Jun Wang.
\newblock Learning to select cuts for efficient mixed-integer programming.
\newblock {\em Pattern Recognition}, 123:108353, 2022.

\bibitem{hussein2017imitation}
Ahmed Hussein, Mohamed~Medhat Gaber, Eyad Elyan, and Chrisina Jayne.
\newblock Imitation learning: A survey of learning methods.
\newblock {\em ACM Computing Surveys (CSUR)}, 50(2):1--35, 2017.

\bibitem{labassi2022learning}
Abdel~Ghani Labassi, Didier Ch{\'e}telat, and Andrea Lodi.
\newblock Learning to compare nodes in branch and bound with graph neural networks.
\newblock {\em Advances in Neural Information Processing Systems}, pages 32000--32010, 2022.

\bibitem{land2010automatic}
Ailsa~H. Land and Alison~G. Doig.
\newblock An automatic method for solving discrete programming problems.
\newblock In Michael J{\"{u}}nger, Thomas~M. Liebling, Denis Naddef, George~L. Nemhauser, William~R. Pulleyblank, Gerhard Reinelt, Giovanni Rinaldi, and Laurence~A. Wolsey, editors, {\em 50 Years of Integer Programming 1958-2008 - From the Early Years to the State-of-the-Art}, pages 105--132. Springer, 2010.

\bibitem{loshchilov2016sgdr}
Ilya Loshchilov and Frank Hutter.
\newblock Sgdr: Stochastic gradient descent with warm restarts.
\newblock {\em arXiv preprint arXiv:1608.03983}, 2016.

\bibitem{loshchilov2017decoupled}
Ilya Loshchilov and Frank Hutter.
\newblock Decoupled weight decay regularization.
\newblock {\em arXiv preprint arXiv:1711.05101}, 2017.

\bibitem{MaherMiltenbergerPedrosoRehfeldtSchwarzSerrano2016}
Stephen Maher, Matthias Miltenberger, Jo{\~{a}}o~Pedro Pedroso, Daniel Rehfeldt, Robert Schwarz, and Felipe Serrano.
\newblock {PySCIPOpt}: Mathematical programming in python with the {SCIP} optimization suite.
\newblock In {\em Mathematical Software {\textendash} {ICMS} 2016}, pages 301--307. Springer International Publishing, 2016.

\bibitem{maros2002computational}
Istv{\'a}n Maros.
\newblock {\em Computational techniques of the simplex method}.
\newblock Springer Science \& Business Media, 2002.

\bibitem{miltenberger2018exploring}
Matthias Miltenberger, Ted Ralphs, and Daniel~E Steffy.
\newblock Exploring the numerics of branch-and-cut for mixed integer linear optimization.
\newblock In {\em Operations Research Proceedings 2017: Selected Papers of the Annual International Conference of the German Operations Research Society (GOR), Freie Universi{\"a}t Berlin, Germany, September 6-8, 2017}, pages 151--157. Springer, 2018.

\bibitem{mnih2016asynchronous}
Volodymyr Mnih, Adria~Puigdomenech Badia, Mehdi Mirza, Alex Graves, Timothy Lillicrap, Tim Harley, David Silver, and Koray Kavukcuoglu.
\newblock Asynchronous methods for deep reinforcement learning.
\newblock In {\em International conference on machine learning}, pages 1928--1937. PMLR, 2016.

\bibitem{nair2020solving}
Vinod Nair, Sergey Bartunov, Felix Gimeno, Ingrid Von~Glehn, Pawel Lichocki, Ivan Lobov, Brendan O'Donoghue, Nicolas Sonnerat, Christian Tjandraatmadja, Pengming Wang, et~al.
\newblock Solving mixed integer programs using neural networks.
\newblock {\em arXiv preprint arXiv:2012.13349}, 2020.

\bibitem{padberg1991branch}
Manfred Padberg and Giovanni Rinaldi.
\newblock A branch-and-cut algorithm for the resolution of large-scale symmetric traveling salesman problems.
\newblock {\em SIAM Review}, 33(1):60--100, 1991.

\bibitem{pateria2021hierarchical}
Shubham Pateria, Budhitama Subagdja, Ah-hwee Tan, and Chai Quek.
\newblock Hierarchical reinforcement learning: A comprehensive survey.
\newblock {\em ACM Computing Surveys (CSUR)}, 54(5):1--35, 2021.

\bibitem{paulus2022learning}
Max~B Paulus, Giulia Zarpellon, Andreas Krause, Laurent Charlin, and Chris Maddison.
\newblock Learning to cut by looking ahead: Cutting plane selection via imitation learning.
\newblock In {\em International Conference on Machine Learning}, pages 17584--17600. PMLR, 2022.

\bibitem{tang2020reinforcement}
Yunhao Tang, Shipra Agrawal, and Yuri Faenza.
\newblock Reinforcement learning for integer programming: Learning to cut.
\newblock In {\em International Conference on Machine Learning}, pages 9367--9376. PMLR, 2020.

\bibitem{wang2023learning}
Zhihai Wang, Xijun Li, Jie Wang, Yufei Kuang, Mingxuan Yuan, Jia Zeng, Yongdong Zhang, and Feng Wu.
\newblock Learning cut selection for mixed-integer linear programming via hierarchical sequence model.
\newblock {\em arXiv preprint arXiv:2302.00244}, 2023.

\bibitem{wesselmann2012implementing}
Franz Wesselmann and Uwe Stuhl.
\newblock Implementing cutting plane management and selection techniques.
\newblock In {\em Technical Report}. University of Paderborn, 2012.

\bibitem{zarpellon2021parameterizing}
Giulia Zarpellon, Jason Jo, Andrea Lodi, and Yoshua Bengio.
\newblock Parameterizing branch-and-bound search trees to learn branching policies.
\newblock In {\em Proceedings of the AAAI Conference on Artificial Intelligence}, volume~35, pages 3931--3939, 2021.

\end{thebibliography}

\newpage
\appendix

\section{Background}

\subsection{Cutting Planes (Cuts)}
Cuts can be used to accelerate the solving process of MILP problems. Given a MILP problem, as shown in Eq.~(1), the solving steps by using the cut method are as follows:
\begin{enumerate}
\item First, ignore the integer constraints of the original MILP problem Eq.~(1) and obtain its LP relaxation:
\begin{equation}        z^{\ast}_{LP}=\min\{\mathbf{c}^T\mathbf{x}|\mathbf{A}\mathbf{x}\le \mathbf{b},\mathbf{x}\in \mathbb{R}^n\},
\end{equation}
where $z^{\ast}_{LP}$ represents the objective value corresponding to the optimal solution $\mathbf{x}^{\ast}_{LP}$ of the current LP relaxation. Let $C_{LP}$ be the feasible region of the current LP relaxation, and $C_{IP}$ be the corresponding feasible region of the original MILP problem. Since the LP relaxation is obtained by ignoring the integer constraints of the original MILP problem, $C_{LP}$ contains $C_{IP}$.
\item Assuming that $\mathbf{x}^{\ast}_{LP}$ does not fully satisfy the integer constraints, we can construct linear inequalities $\alpha x\leq \beta$ to tighten the feasible region without deleting any integer feasible solutions of the original problem Eq.~(1). These linear inequalities are called cuts, which can be added to the problem Eq.~(3) without loss of optimal solutions.
\item Solve Eq.~(3) with new constraints, namely $\mathbf{A}\mathbf{x}\le \mathbf{b}$ and $\alpha x\leq \beta$. If the optimal value $z^{\ast}_{LP}$ is reached, the solving process ends. If not,  repeatedly execute Step 2 to tighten the feasible region and eventually reach the optimal solution of the MILP problem at an integer extreme point within the narrowed feasible region.
\end{enumerate}
In modern MILP solvers, different methods are used to generate a set of candidate cuts by adjusting the solver's parameters~\cite{dey2018theoretical}. However, adding too many cuts may increase the problem complexity, leading to computational burden and numerical instability~\cite{wesselmann2012implementing,dey2018theoretical}. Therefore, selecting appropriate cuts from the candidate set is important for the efficiency of MILP solving.

\subsection{Branch and Cut}
Due to numerical instability, the cut method may be ineffective if it is used solo. Therefore, advanced solvers typically integrate the cut method into B\&B to form B\&C~\cite{miltenberger2018exploring}. To solve problem $P$, B\&C starts by gradually constructing a B\&B search tree $T$ from the root node. In the initial stage, B\&B finds a feasible solution $\mathbf{x}_h$ based on a specific heuristic algorithm, where $\mathbf{x}_h$ satisfies the constraints in Eq.~(1) but may not necessarily minimize $f(\mathbf{x})=\mathbf{c}^T\mathbf{x}$. Let $B=f(\mathbf{x}_h)$ denote the upper bound of the problem, where $B$ represents the objective value corresponding to the best solution found so far. Subsequently, a subproblem is selected from the list of unexplored subproblems. If a feasible solution $\mathbf{x}_h^{\prime}$ is found in the solution space of this subproblem such that $f(\mathbf{x}_h^{\prime})\textless f(\mathbf{x}_h)$, then $B$ is updated to $B=f(\mathbf{x}_h^{\prime})$, and $\mathbf{x}_h^{\prime}$ becomes the best solution found so far. If it can be proven that there is no feasible solution $\mathbf{x}_h^{\prime}$ in the solution space of this subproblem such that $f(\mathbf{x}_h^{\prime})\textless f(\mathbf{x}_h)$, then this subproblem is pruned. Otherwise, the subproblem is further divided, and the node representing the subproblem is added to the search tree $T$. This process is repeated until there are no unexplored subproblems. At this point, $\mathbf{x}_h$ is the optimal solution, and $B$ represents the objective function value corresponding to the optimal solution. 

\section{Related work}
In recent years, 
the integration of deep learning techniques into traditional solving methods for MILP problems has become a new trend~\cite{gasse2019exact,gupta2020hybrid,zarpellon2021parameterizing,nair2020solving,labassi2022learning,paulus2022learning,tang2020reinforcement,wang2023learning,huang2022learning}. Among these approaches, techniques such as imitation learning~\cite{hussein2017imitation} are widely used to improve decision-making in B\&B process. While graph neural networks have previously been used to solve specific combinatorial optimization problems, Gasse et al.~\cite{gasse2019exact} were the first to introduce them into the variable selection process. They opted for graph neural networks to better exploit the bipartite structure of MILP problems by employing heuristic algorithms as expert policy and imitating it through behavioral cloning, which is simpler than predicting variable scores or rankings. Labassi et al.~\cite{labassi2022learning} followed a similar approach and applied the same idea and model to the node selection process.

For cut selection, Paulus et al.~\cite{paulus2022learning} proposed a forward-looking heuristic cut selection strategy and employed the same behavioral cloning approach and a similar graph neural network to mimic this strategy. The mentioned works are all based on the idea of imitating high-quality expert strategies, while Bengio et al.~\cite{bengio2021machine} argued that imitation learning alone could be valuable if the imitated strategy is as good as the quality of the expert strategy and significantly accelerates computation. Conversely, if given enough training, models trained based on reward signals may surpass expert strategies. Tang et al.~\cite{tang2020reinforcement} modeled the cut selection problem as a Markov decision process and designed a deep reinforcement learning architecture to explore more efficient strategies. Similarly, Wang et al.~\cite{wang2023learning} employed reinforcement learning methods to explore new strategies. They noticed that both the number and order of selected cuts can impact the quality of the solution, thus transforming the cut selection problem into a sequence-to-sequence problem. Additionally, Huang et al.~\cite{huang2022learning} utilized a multi-instance learning approach to measure the ranking of cuts and make selections based on that.

It is remarkable that the method proposed in this paper is inspired by the aforementioned ideas of leveraging graph neural networks to capture the underlying structure of MILP problems and use reinforcement learning methods to explore new selection strategies. 
Additionally, the problem of cut order is neglected in previous works. Since it has big impact on MILP solving, in this paper, we solve this problem and make our model stable whatever the order of cuts. 

\section{Motivation Supplement}
\label{MoreMotivation}
\begin{figure}[H]
\centering
\subfigure[ ]{
\includegraphics[width=0.48\textwidth]{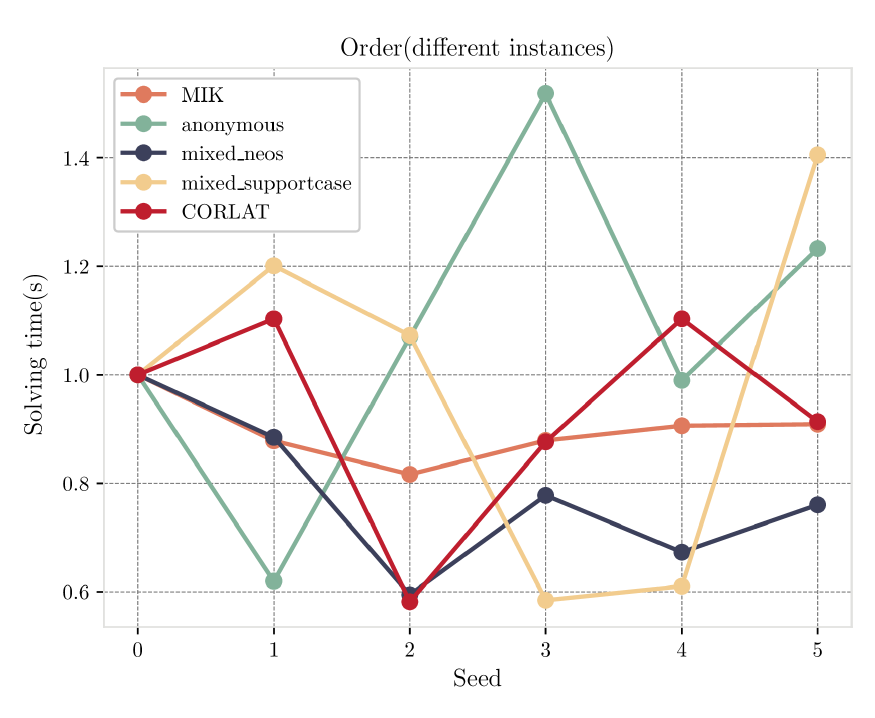}}
\subfigure[ ]{
\includegraphics[width=0.48\textwidth]{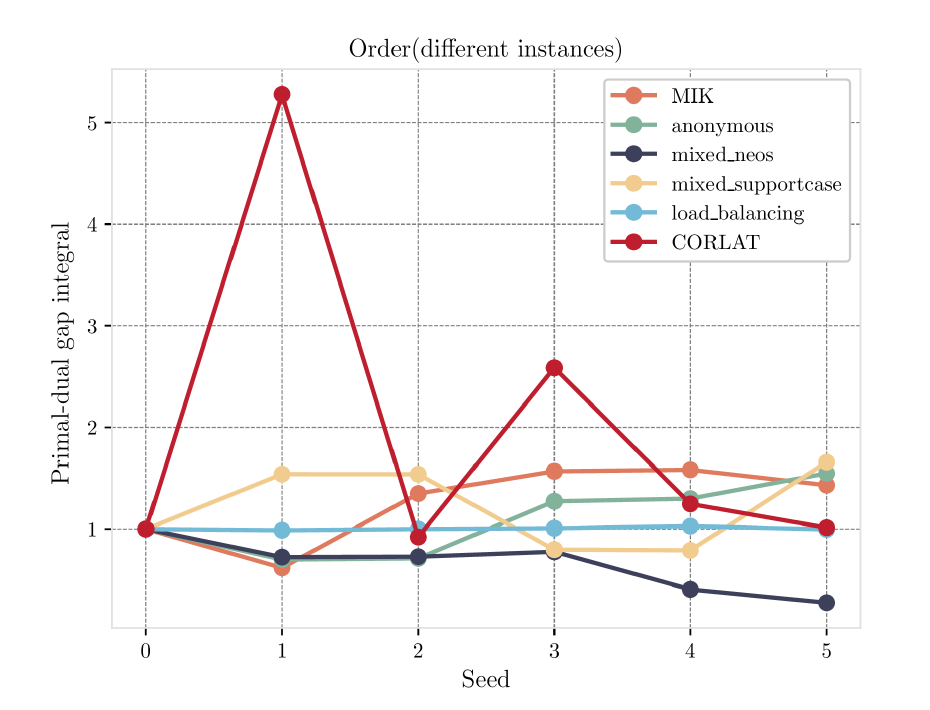}}
\caption{(a) From five datasets, take a MILP problem from each dataset and measure the solving time when the input cut sequences are shuffled by using different random seeds. (b) Primal-dual gap integral was obtained using the same method. In both figures, the horizontal ordinate $0$ represents the use of the default input sequence from SCIP without modifications.}
\label{motivation-more-ins}
\end{figure}

\begin{figure}[H]
\centering
\subfigure[ ]{
\includegraphics[width=0.48\textwidth]{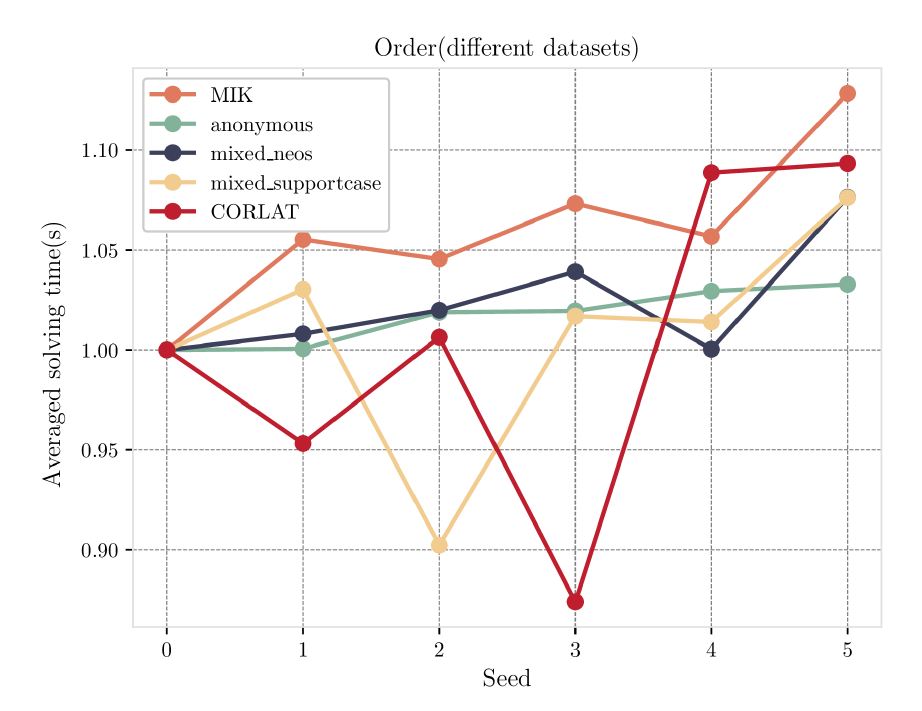}}
\subfigure[ ]{
\includegraphics[width=0.48\textwidth]{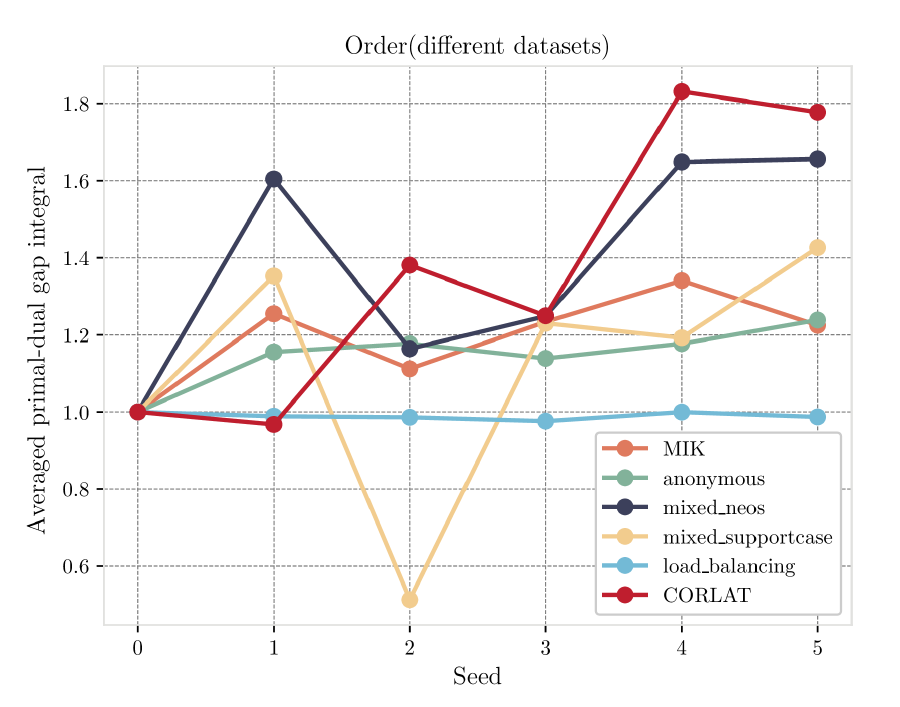}}
\caption{(a) Average solving time obtained by shuffling the input cut sequences with different random seeds for each problem from the five datasets. (b) Primal-dual gap integral was obtained using the same method. In both figures, the horizontal ordinate $0$ represents the use of the default input sequence from SCIP without modifications.}
\label{motivation-more-set}
\end{figure}

The experimental results are shown in Figure~\ref{motivation-more-ins} and Figure~\ref{motivation-more-set}. From the experimental results, one can observe that: (1) For different types of MILP problems, changing the cut order will significantly and irregularly affect the problem-solving process. (2) For different MILP datasets, although increasing the number of problems may mitigate the impact of cut sequence changes, it still results in a certain degree of variation in the solution quality. Therefore, the experimental results demonstrate that the sequence-based selection strategy captures positional information of the input cuts, which interferes with problem solving and reduces solution stability.

\begin{table}
\caption{The sequence of cuts is generated by the sequence model and then shuffled with different random number seeds. Here, ``default'' represents the sequence of cuts generated by SCIP without any modification.}
\label{motivation}
\centering
\scalebox{0.8}{
\begin{tabular}{ll}
\toprule
Seed & Selected sequence of cuts \\
\midrule
\multirow{2}*{default} & 0, 8, 48, 41, 43, 49, 44, 1, 9, 6, 2, 5, 3, 4, 47, 7, 45, 46, 11, 12, 42, 10, 50, 40, 17, 13, 32 \\ 
& 16, 15, 31, 24, 22, 35, 14, 26, 25, 37, 36, 21, 23, 20, 18, 38 \\
\midrule
\multirow{2}*{1} & 0, 8, 2, 1, 3, 43, 5, 49, 41, 4, 44, 6, 48, 7, 9, 15, 47, 40, 27, 46, 45, 11, 13, 24, 16, 36, 21 \\ 
& 17, 14, 30, 31, 32, 12, 25, 26, 18, 22, 35, 37, 42, 28, 50, 23 \\
\midrule
\multirow{2}*{2} & 0, 8, 1, 15, 2, 44, 5, 49, 4, 48, 43, 3, 6, 41, 7, 9, 36, 33, 29, 47, 45, 40, 35, 12, 17, 21, 11 \\
& 34, 46, 32, 14, 50, 20, 24, 25, 13, 26, 18, 22, 37, 31, 16, 42 \\
\midrule
\multirow{2}*{3} & 8, 0, 26, 4, 1, 6, 44, 5, 43, 48, 3, 2, 49, 41, 7, 45, 50, 9, 34, 40, 11, 47, 37, 12, 17, 24, 46\\
& 15, 42, 21, 32, 14, 18, 31, 35, 30, 16, 27, 25, 22, 39, 36, 38 \\
\midrule
\multirow{2}*{4} & 0, 8, 48, 36, 6, 2, 3, 1, 44, 49, 43, 41, 4, 5, 7, 37, 22, 9, 26, 47, 11, 17, 40, 46, 38, 42, 31\\ 
& 24, 45, 15, 19, 21, 27, 13, 12, 25, 16, 32, 35, 18, 30, 23, 20 \\
\midrule
\multirow{2}*{5} & 0, 8, 6, 11, 5, 4, 1, 43, 49, 48, 41, 2, 3, 44, 7, 9, 35, 13, 18, 47, 46, 24, 45, 40, 22, 37, 25\\ 
& 17, 16, 12, 26, 36, 15, 32, 50, 38, 31, 42, 21, 33, 10, 30, 39 \\
\bottomrule
\end{tabular}
}
\end{table}

\section{Methodology Supplement}
\subsection{Feature Details}
\label{FeatureDetails}
We use SCIP as the solver for MILP problems and collect state information during the problem-solving process by using custom PySCIPOpt functions. The state information consists of three parts: variable features, constraint features, and cut features. 
\begin{enumerate}
    \item For variable features, we follow the setting in \cite{paulus2022learning}, and remove two features related to incumbent values, like the setting in  \cite{gasse2019exact}.
    \item We notice that some constraint features designed in \cite{paulus2022learning} are more relevant for cuts. Hence, we exclude these features in constraint features.
    \item We adhere to the designs of cut features in \cite{huang2022learning, wesselmann2012implementing, dey2018theoretical, wang2023learning} without modification.
\end{enumerate}
Note that we did not manually design features for the edges of associated nodes. Conversely, we treat the edge values as attention scores between adjacent nodes, which are calculated through HGT.
We provide detailed descriptions of variable features, constraint features, and cut features in Table~\ref{feature}.

\begin{table}
\caption{The features of variables, constraints,
and cuts.}
\label{feature}
\centering
\scalebox{0.8}{
\begin{tabular}{lll}
\toprule
Node&Feature&Description \\
\hline
\multirow{11}*{\textbf{Vars}} & norm\_coef & Objective coefficient, normalized by objective norm \\
~    & type     & Variable type (binary, integer, implint, continuous), one-hot\\
~    & has\_lb    & Lower bound in LP\\
~    & has\_ub    & Upper bound in LP\\
~    & norm\_redcost    & Reduced cost value in LP, normalized by objective norm\\[0.5ex]
~    & solval    & LP solution of a variable\\
~    & solfrac   & Fractional part of value\\
~    & sol\_is\_at\_lb    & Solution value equals lower bound\\
~    & sol\_is\_at\_ub    & Solution value equals upper bound\\
~    & norm\_age  & Number of successive times, normalized by number of solved LPs\\
~    & basestat  & Basis status (lower, basic, upper, zero), one-hot\\
\midrule
\multirow{13}*{\textbf{Cons}} & rank & Rank of a row \\
~ & norm\_nnzrs & Nonzero entries in row vector \\
~ & norm\_bias & Unshifted side, normalized by row norm \\
~ & row\_is\_at\_lhs & Row value equals left hand side \\
~ & row\_is\_at\_rhs & Row value equals right hand side \\
~ & norm\_dualsol & Dual LP solution of a row, normalized by row and objective norm \\
~ & basestat &  Basis status (lower, basic, upper, zero), one-hot \\
~ & norm\_age & Number of successive times, normalized by total number of solved LPs \\
~ & norm\_nlp\_creation & LPs since the row has been created, normalized \\
~ & norm\_intcols & Number of integral columns in the row \\
~ & is\_integral & Activity of the row is always integral in a feasible solution \\
~ & is\_removable & Row is removable from the LP \\
~ & is\_in\_lp & Row is member of current LP \\
\midrule
\multirow{7}*{\textbf{Cuts}} & cut\_coefficients & Mean, max, min, std of cut coefficients \\
~ & objective\_coefficients & Mean, max, min, std of the objective coefficients \\
~ & parallelism & Parallelism between the objective and the cut $\frac{c^T\alpha}{\left|c\right|\left|\alpha\right|}$ \\
~ & efficacy & Euclidean distance of the cut hyperplane to the current LP solution \\
~ & support & Proportion of non-zero coefficients of the cut\\
~ & integral\_support & Proportion of non-zero coefficients with respect to integer variables of the cut \\
~ & normalized\_violation & Violation of the cut to the current LP solution $\max\{0,\frac{\alpha^Tx^*_{LP}-\beta}{\left|\beta\right|}\}$ \\
\bottomrule
\end{tabular}
}
\end{table}

\subsection{HGT Details}
\label{HGT}
\begin{figure}
\centering
\includegraphics[width=0.95\textwidth]{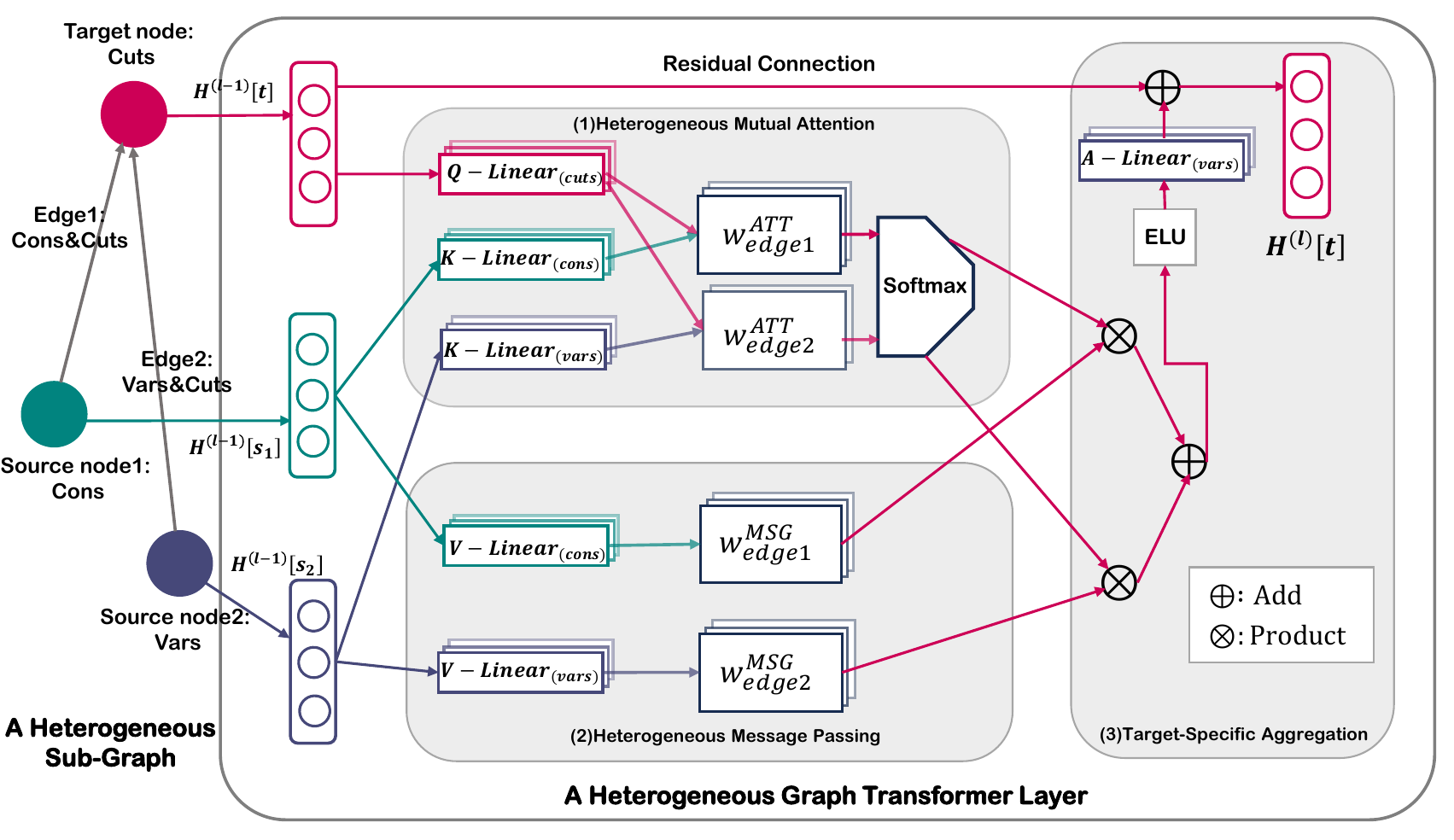}
\caption{\footnotesize A demo of HGT extracting MILP state information from the tripartite graph. The heterogeneous sub-graph is a directed graph from variables and constraints to cuts. HGT learns the context representation of each node based on different meta-relations as input, which is used for subsequent cut selection tasks.}
\label{hgt}
\end{figure}

Drawing inspiration from the classic Transformer, the Heterogeneous Graph Transformer(HGT) encodes target nodes as query vectors and source nodes as key vectors and value vectors. By computing the attention scores of source nodes to target nodes, corresponding weights are obtained for message propagation and aggregation, updating the feature vectors of each node. A single-layer HGT network structure is illustrated in Figure~\ref{hgt}. It specifically treats different types of edges and nodes in different networks, which aligns perfectly with our encoding of MILP states as heterogeneous tripartite graphs. We primarily leverage HGT to extract underlying information from the MILP problem for subsequent cut selection tasks. 

\section{Experiment Supplement}

\subsection{Datasets Details}
\label{DatasetsDetails}
For convenient comparison, the datasets of the MILP problem are the same as those in ~\cite{wang2023learning}. However, we deviate by excluding the easy dataset of simple instances (problems that can be solved optimally by SCIP 8.0.0 within one minute) and focus on more challenging MILP problems (medium dataset and hard dataset). These datasets are described as follows. 

\textbf{Medium datasets.} The SCIP 8.0.0 solver needs at least five minutes to get the optimal value for a MILP problem in the medium dataset.
\begin{enumerate}
    \item MIK.
    Similar to the works \cite{he2014learning, wang2023learning}, we include MIK as one of the problem sets. MIK comprises a set of MILP problems with knapsack constraints, and further details about MIK can be found at https://atamturk.ieor.berkeley.edu/data/mixed.integer.knapsack/.
    \item CORLAT.
    Similar to the works \cite{he2014learning, nair2020solving, wang2023learning}, we include CORLAT as one of the problem sets. CORLAT is a real-world, small-scale problem set used for wildlife management. The number of variables ranges from 300 to 1000, and the number of constraints ranges from 100 to 500. More details about CORLAT can be found at https://bitbucket.org/mlindauer/aclib2/src/master/.
\end{enumerate}
\textbf{Hard Datasets.} For a MILP problem in the hard dataset, the SCIP 8.0.0 solver needs at least one hour to obtain the optimal value.
\begin{enumerate}
    \item 
     MIPLIB mixed neos and MIPLIB mixed supportcase. Two problem sets are subsets of MIPLIB 2017. They are collections of problem sets with similar structures constructed based on manually designed features by \cite{wang2023learning}.
    \item 
    Load Balancing and Anonymous. Two problem sets are from the NeurIPS 2021 ML4CO competition.  More details about the problem sets can be obtained from https://www.ecole.ai/2021/ml4co-competition/.
\end{enumerate}

Each dataset is divided into a training set comprising $80\%$ of the problem instances and a testing set comprising $20\%$ of the instances. 
Table~\ref{datasets} provides detailed information on the scale of different problem sets. All datasets in this paper can be downloaded from the anonymous webpage, \url{https://figshare.com/s/0b1f463c37e3f285439d}.

\begin{table}
\caption{The average number of constraints and variables in different datasets, where $Cons$ represents constraints and $Vars$ represents variables.}
\label{datasets}
\centering
\scalebox{0.8}{
\begin{tabular}{lllllll}
\toprule
Datasets & MIK & CORLAT & MIPLIB mixed neos & MIPLIB mixed supportcase & Load Balancing & Anonymous \\
\midrule
$Cons$ & 346 & 486 & 5,660 & 19,910 & 64,304 & 49,603 \\
$Vars$ & 413 & 466 & 6,958 & 19,766 & 61,000 & 37,881  \\
\bottomrule
\end{tabular}
}
\end{table}

\subsection{Baselines Details}
\label{Baselines}
In this section, we provide a detailed description of all baselines, including heuristic methods and learning-based methods. Table~\ref{baseline} presents heuristic cut selection strategies used as baselines. All heuristic strategies are implemented through PySCIPOpt and SCIP. Except for the default SCIP strategy, the proportion of cuts selected by other heuristic strategies is fixed at $0.2$, following the setting in \cite{wang2023learning}. We also consider the SOTA cut selection method HEM proposed by Wang et al.~\cite{wang2023learning} as a baseline. We conduct comparative experiments by reproducing their open-source code with same hyperparameters set as the original paper.

\begin{table}
\caption{Detailed description of different baselines.}
\label{baseline}
\centering
\scalebox{0.8}{
\begin{tabular}{lll}
\toprule
Method & Description \\
\midrule
No Cuts & Do not add any cuts. \\
\midrule
Random &  Random selects a fixed ratio of the candidate cuts stochastically. \\
\midrule
\multirow{3}*{Normalized Violation} & Normalized Violation is a score-based rule. It scores each cut based on the normalized \\
 & violation of the cuts to the current LP solution, and selects a fixed ratio of cuts \\
 & with high scores. The normalized violation is defined by $\max\{0,\frac{\alpha^Tx^*_{LP}-\beta}{\left| \beta\right|}\}$.\\
\midrule
\multirow{2}*{Efficacy} & Efficacy is a score-based rule. It scores each cut based on the Euclidean distance of the cut \\
 & hyperplane to the current LP solution, and selects a fixed ratio of cuts with high scores.\\
\midrule
\multirow{2}*{Default} & Default is the default cut selection rule used in SCIP 8.0.0 (Bestuzheva et al.~\cite{bestuzheva2021scip}). Default \\
 & selects variable ratios of cuts rather than a fixed ratio. \\
\midrule
HEM & The SOTA sequence-model-based cuts selection method provided by ~\cite{wang2023learning}.\\
\bottomrule
\end{tabular}
}
\end{table}

\subsection{Training Algorithm}
\label{TrainAlgo}
We present the training process of our model, as shown in Algorithm ~\ref{training}.
\begin{algorithm}
	\renewcommand{\algorithmicrequire}{\textbf{Input:}}
	\renewcommand{\algorithmicensure}{\textbf{Output:}}
	\caption{The training process of our model.}
	\label{training}
	\begin{algorithmic}[1]
		\REQUIRE MILP problem set $D_m$, training epoch $N_e$, batch size $N_b$, learning rate $\alpha$.
        \STATE Initialize policy network parameters $(\theta_h,\theta_l)$, higher-layer network training datasets $D^h_{train}$, lower-layer network training datasets $D^l_{train}$;
        \FOR{$N_e$ epochs}
        \STATE Empty training datasets $D^h_{train}$ and $D^l_{train}$;
        \FOR{$N_b$ steps}
        \STATE Solve MILP instances, collect status information $s_0$ from $D_m$;
        \STATE Construct a tripartite graph $G$ as network input;
        \STATE Take action ${ratio, a_0}$ based on the graph $G$ and policy $\pi$;
        \STATE Obtain reward $r_0$ and construct $(s_0, ratio, r_0)$ to $D^h_{train}$ and $(s_0, a_0, r_0)$ to $D^l_{train}$; 
        \ENDFOR
        \STATE Compute policy gradient using $D^h_{train}$ and $D^l_{train}$;
        \STATE Update parameters, $\theta_h=\theta_h+\alpha\nabla_{\theta_h}J([\theta_h,\theta_l]), \theta_l=\theta_l+\alpha\nabla_{\theta_l}J([\theta_h,\theta_l])$;
        \ENDFOR
	\end{algorithmic}  
\end{algorithm}

\subsection{Experiment Setting}
Adding cuts to the original MILP problem before the branching process is crucial for enhancing the relaxation of linear programming and improving solving efficiency~\cite{wesselmann2012implementing,bengio2021machine}. Similar to the works ~\cite{gasse2019exact, huang2022learning, paulus2022learning, wang2023learning}, we allow cuts to be added at the root node and set the number of rounds of adding cuts to $1$.  During this process of MILP solving, we can collect training data $(s_0, a_0, r_0)$, which implies that after Monte Carlo approximation, the action-value function $Q_{\pi(\theta)}(s_0,a_0) = r(s_0,a_0)$. Additionally, we use a neural network to approximate the value function $V_{\pi(\theta)}$ to compute the advantage function. Each trajectory generated by MILP has a length of 1, and we train with batches of $16$ data points. We present the algorithmic procedure for training in Appendix \ref{TrainAlgo}.

All experiments are conducted using  SCIP 8.0.0 as the solver for MILP problems. We control the process of data collection and problem-solving using the PySCIPOpt 4.1.0 interface library~\cite{MaherMiltenbergerPedrosoRehfeldtSchwarzSerrano2016}, which integrates SCIP with Python. We utilize 17 different types of separators available in SCIP, including \textit{aggregation, clique, close-cuts, CMIR, Convex Projection, disjunctive, edge-concave, flow cover, gauge, Gomory, implied bounds, integer objective value, knapsack cover, multi-commodity-flow network, odd cycle, rapid learning, and zero-half}.

Since there are no existing functions available in the PySCIPOpt library for feature extraction, we implemented two functions, \textit{getVarConState} and \textit{getVarCutCoef}, to help extract the required features. We follow the experimental settings in \cite{gasse2019exact, huang2022learning, paulus2022learning, wang2023learning}, which allow only the addition of cuts at the root node and set the number of rounds for adding cuts to $1$. For other parameters, we keep the default settings in SCIP to ensure a fair experimental environment, and we enable SCIP's presolving and heuristics functionalities to assist in solving while ensuring consistency between the experimental and real environments.
For all experiments, we set a maximum solving time limit of $300$ seconds. Additionally, since there are problems in all datasets that cannot be solved within the specified time, we set the reward function for all problems to be the primal-dual gap integral.
We use the AdamW optimizer~\cite{loshchilov2017decoupled} with a Cosine Annealing Learning Rate Scheduler~\cite{loshchilov2016sgdr} for optimization. The remaining hyperparameters follow the settings of HEM~\cite{wang2023learning}.

\subsection{Hardware Specifics}
All experiments are conducted on a machine equipped with an RTX 4090Ti GPU, an Intel i9-13900KF CPU and 64GB memory.

\subsection{Codes}
The implementation codes of our model HGTSM can be downloaded from the anonymous webpage, \url{https://figshare.com/s/0b1f463c37e3f285439d}. 


\section{Experiments and Results Supplement}

\subsection{More Stability Result}
\label{MoreStaResult}
\begin{figure}[H]
\centering
\subfigure[ ]{
\includegraphics[width=0.48\textwidth]{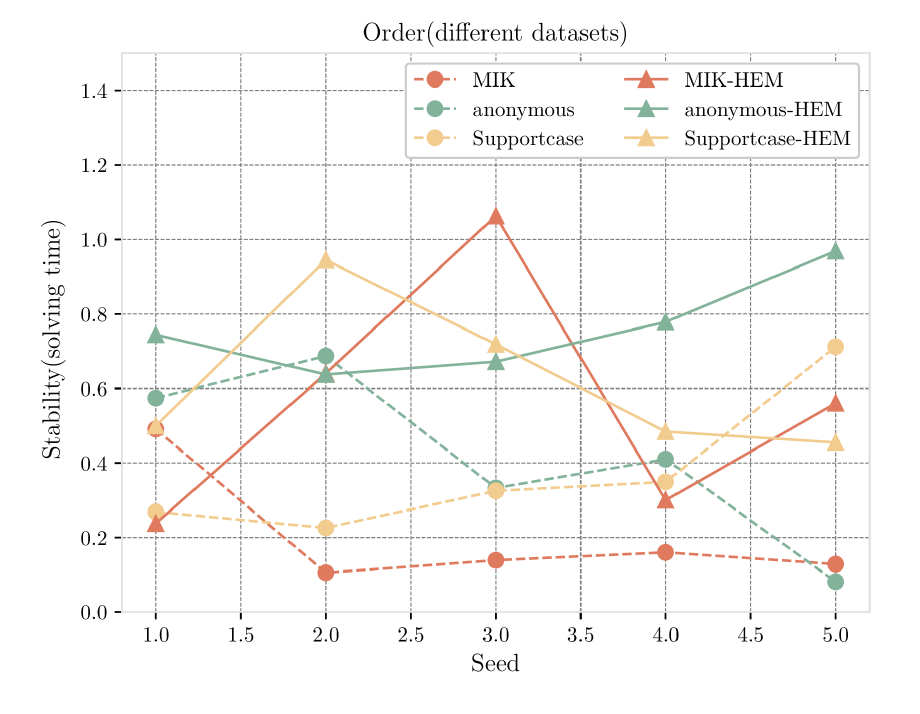}}
\subfigure[ ]{
\includegraphics[width=0.48\textwidth]{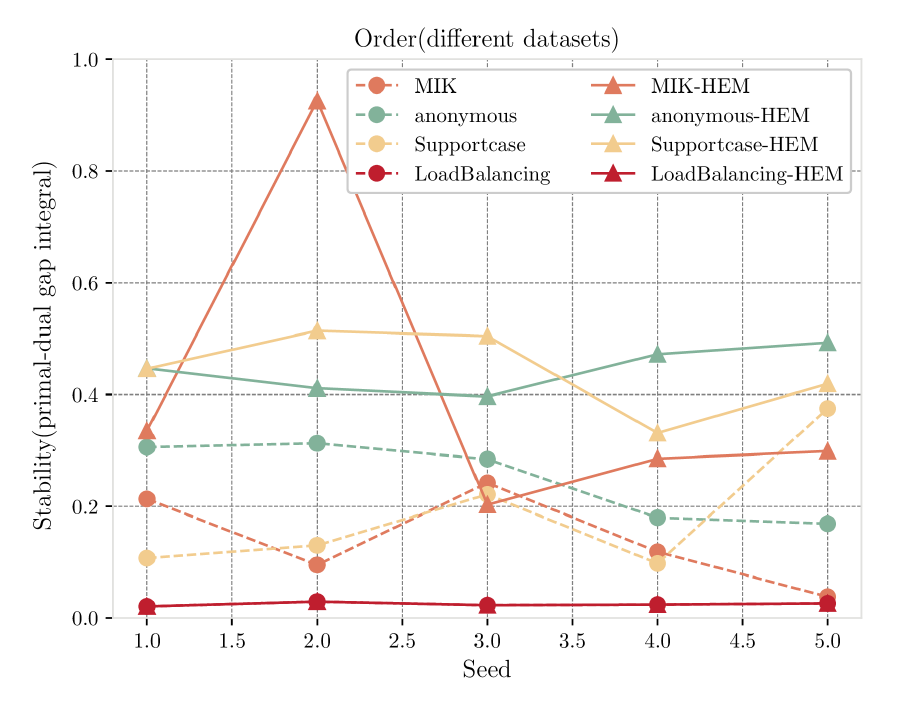}}
\caption{(a) Comparison of Stability (solving time) between HEM and HGTSM across different datasets under various random seeds. (b) Comparison of Stability (primal-dual gap integral) between HEM and HGTSM across different datasets under various random seeds.  In both figures, the same color represents the same problem, with HEM represented by straight lines and HGTSM represented by line segments.}
\label{sta-more-set}
\end{figure}

We separately evaluate the trends in stability across different datasets under various random seeds. Due to the limited number of instances in the MIPLIB mixed Neos test set and all instances in the Load Balancing test set exceeding the set solving time limit, we report the trends in solving time stability for three additional datasets and the trends in primal-dual gap integral for four datasets (CORLAT has been provided in the Section~\ref{ExpSta}). The experimental results are illustrated in the Figure~\ref{sta-more-set}. From the experimental results, it is evident that except for the two models performing equally well on the Load Balancing dataset, our model exhibits less fluctuation in solving across other datasets, with lower $Stability$ scores. This indicates a more pronounced advantage of stability for our model compared to HEM~\cite{wang2023learning}.

\subsection{More Generalization Result}
\label{MoreGenResult}
In this section, we report more results on the generalization experiment. Besides the test of primal-dual gap integral shown in Table~\ref{result_Gen}, we use models trained on a specific dataset to test other datasets of MILP problems and record the mean and standard deviation of the solving time. The experimental results are shown in Table~\ref{result_gen_t}. Surprisingly, and in contrast to the situation with the result of primal-dual gap integral (Table~\ref{result_Gen}), the original models do not always win in terms of solving time for every problem. In the datasets of CORLAT, Anonymous, MIPLIB mixed neos, and MIPLIB mixed supportcase, models trained on other problems outperformed the original models in solving time. On the other hand, the reward function for all problems was set to the primal-dual gap integral, which might be one of the reasons for this outcome. Considering the primal-dual gap integral, the models for specific tasks still have a significant advantage. Nevertheless, the experimental results still motivate us to seek more efficient and general methods for solving MILP problems.

\begin{table}
\caption{The statistical comparison of four datasets includes the mean (standard deviation) of solving time. Note that Sta. means stability, and Time Imp. means the improvement(\%) of solving time.}
\label{result_2_2}
\centering
\scalebox{0.74}{
\begin{tabular}{lllllllll}
\toprule
~ & \multicolumn{2}{c}{CORLAT} & \multicolumn{2}{c}{MIK} & \multicolumn{2}{c}{Anonymous} & \multicolumn{2}{c}{MIPLIB mixed upporecase}\\
\cmidrule(r){1-3} \cmidrule(r){4-5} \cmidrule(r){6-7} \cmidrule(r){8-9}
Method  & Sta. $\downarrow$ & Time Imp. $\uparrow$ & Sta. $\downarrow$ & Time Imp. $\uparrow$ & Sta. $\downarrow$ & Time Imp. $\uparrow$ & Sta. $\downarrow$ & Time Imp. $\uparrow$  \\
\cmidrule(r){1-3} \cmidrule(r){4-5} \cmidrule(r){6-7} \cmidrule(r){8-9}
HEM & 1.35 & NA & 0.56 & NA & 0.76 & NA & 0.62 & NA \\
HGTSM & \textbf{0.15} & \textbf{88.45} & \textbf{0.21} & \textbf{63.40} & \textbf{0.41} & \textbf{45.11} & \textbf{0.37} & \textbf{39.34} \\
\bottomrule
\end{tabular}
}
\end{table}

\begin{table}
\caption{The generalization results of our model on the metric of time. For example, $53.06(105.99)$ in the cell $[2,1]$ of datagrid means that the model is trained on the dataset MIK, and evaluated on the dataset CORLAT. $53.06$ is the average of solving times, and $105.99$ is the standard deviation of solving times. Note that Neos means MIPLIB mixed neos, and Supportcase means MIPLIB mixed supportcase.}
\label{result_gen_t}
\centering
\scalebox{0.75}{
\begin{tabular}{lllllll}
\toprule
~ & \multicolumn{6}{c}{Solving time using different problem models} \\
\cmidrule(r){2-7}
~ & CORLAT & MIK & Anonymous & Neos & Supportcase & Load Balancing\\
\midrule
CORLAT & 56.23(107.51) & 218.80(107.93) & 238.52(106.76) & \textbf{241.05(102.10)} & 160.21(139.85) & 300.01(0.02)\\
MIK & \textbf{53.06(105.99)} & \textbf{156.60(129.93)} & 243.80(97.67) & 251.16(84.59) & 160.15(140.03) & 300.02(0.02) \\
Anonymous & 58.06(114.42) & 203.12(118.93) & 241.07(104.20) & 250.00(86.59) & 155.20(121.35) &  300.06(0.09)\\
Neos & 84.53(128.40) & 216.25(111.74) & \textbf{233.02(113.60)} & 246.44(89.29) &  \textbf{139.87(130.74)} & 300.02(0.04)\\
Supportcase & 55.96(109.54) & 243.93(95.93) & 250.09(88.28) & 266.56(57.91) & 164.46(136.77) & 300.05(0.12) \\
Load Balancing & 59.31(108.98) & 196.27(113.46) & 238.47(106.96) & 251.32(84.30) & 173.91(131.14) & \textbf{300.01(0.02)} \\
\bottomrule
\end{tabular}
}
\end{table}

\end{document}